\documentclass[12pt]{article}
\usepackage[T1]{fontenc}
\usepackage[dvips]{graphicx}
\usepackage{makeidx}\makeindex
\usepackage[round]{natbib}
\graphicspath{{images/}}
\usepackage[utf8]{inputenc}

\usepackage{amsmath}
\usepackage{geometry}
\usepackage{bbm}
\usepackage{authblk}
\usepackage{amsfonts}
\usepackage{xcolor}
\usepackage{tikz}
\usepackage{comment}
\usepackage{amsthm}
\usepackage{appendix}
\usepackage{bm}
\usepackage{bbm}
\usepackage{algorithm}

\usepackage{mathtools}
\usepackage{amssymb}
\usepackage{enumitem}
\usepackage{verbatim}
\usepackage{amsmath,amsthm}
\usepackage{xspace}
\usepackage{pifont}
\usepackage{graphicx}
\usepackage{amssymb}
\usepackage{epic, eepic}
\usepackage{dsfont}
\usepackage{amssymb}
 \usepackage{mathrsfs}
\usepackage{exscale}
\usepackage{color} 
\usepackage{overpic} 
\usepackage{bm}
\usepackage{bbm}
\usepackage{booktabs} % To thicken table lines
\usepackage{color, colortbl}
\usepackage{subcaption}
\usepackage[round]{natbib}
\usepackage{tikz}
\usepackage{pgfplots}

\usepgfplotslibrary{fillbetween}
\usetikzlibrary{automata,arrows,positioning,calc}
\RequirePackage[colorlinks,citecolor=blue,urlcolor=blue]{hyperref}

\definecolor{Gray}{gray}{0.9}

\usepackage{amsmath,afterpage}
\usepackage{epsf}
\usepackage{graphics,color}

\def\eps{\varepsilon}

\def\lam{\lambda}
\def\rr{\rightarrow}

\def \< {\langle}
\def \> {\rangle}

\def\ol{\overline}

\def\beqa{\begin{eqnarray}}
\def\eeqa{\end{eqnarray}}
\def\beqas{\begin{eqnarray*}}
\def\eeqas{\end{eqnarray*}}

\DeclareMathOperator*{\argmax}{arg\,max}

\newtheorem{theorem}{Theorem}[section]
\newtheorem{lemma}[theorem]{Lemma}

\newtheorem{proposition}[theorem]{Proposition}

\newtheorem{remark}[theorem]{Remark}

\newtheorem{definition}[theorem]{Definition}
\newtheorem{assumption}[theorem]{Assumption}
\numberwithin{equation}{section}
\newcommand{\hatd}[1]{{}}

\newcommand{\sfD}{\mathsf{D}}
\newcommand{\sfE}{\mathsf{E}}

\newcommand{\N}{{\mathbb N}}

\newcommand{\R}{{\mathbb R}}

\newcommand{\indicator}[1]{ \mathbbm{1}_{\{#1\}} }
\newcommand{\calA}{{\cal A}}

\newcommand{\calC}{{\cal C}}
\newcommand{\calD}{{\cal D}}

\newcommand{\calK}{{\cal K}}

\newcommand{\calQ}{{\cal Q}}
\newcommand{\calR}{{\cal R}}
\newcommand{\calS}{{\cal S}}

\newcommand{\calV}{{\cal V}}
\newcommand{\calW}{{\cal W}}

\newcommand{\Backslash}{\text{\textbackslash}}
\newcommand{\PP}{{\mathbb P}}
\newcommand{\EE}{{\mathbb E}}

\newcommand{\bd}{\begin{displaymath}}
\newcommand{\ed}{\end{displaymath}}
\newcommand{\be}{\begin{equation}}
\newcommand{\ee}{\end{equation}}
\newcommand{\bq}{\begin{eqnarray}}
\newcommand{\eq}{\end{eqnarray}}
\newcommand{\bn}{\begin{eqnarray*}}
\newcommand{\en}{\end{eqnarray*}}
\newcommand{\dl}{\delta}

\def\wt{\widetilde}

\usetikzlibrary{decorations.pathreplacing,positioning, arrows.meta}

\usepackage{authblk}

\author[]{Johannes Muhle-Karbe}
\author[]{Eyal Neuman}
\author[]{Yonatan Shadmi}

\affil[]{Department of Mathematics, Imperial College London}
\title{Fluid-Limits of Fragmented Limit-Order Markets}
\date{\today}

\pgfplotsset{compat=1.14} 
\defcitealias{Maglaras2021}{MMZ} % in brackets \citepalias{Maglaras2021}{MMZ21}

\begin{document}
\maketitle
\begin{abstract}
\cite*{Maglaras2021} have introduced a flexible queueing model for fragmented limit-order markets, whose fluid limit remains remarkably tractable. In the present study we prove that, in the limit of small and frequent orders, the discrete system indeed converges to the fluid limit, which is characterized by a system of coupled nonlinear ODEs with singular coefficients at the origin. Moreover, we establish that the fluid system is asymptotically stable for an arbitrary number of limit order books in that, over time, it converges to the stationary equilibrium state studied by~\cite{Maglaras2021}.

\end{abstract}
\begin{description}
\item[Mathematics Subject Classification (2010):] 	60K25, 90B22, 37A50
\item[JEL Classification:] C02, C62, C68 
\item[Keywords:] order routing, limit order books, fragmented markets, fluid limits, asymptotic stability 

\end{description}

\section{Introduction}
Over the last two decades, electronic equity markets have become more and more fragmented rather than gradually converging to a single centralized limit-order book. One important driver for this was the SEC’s ``Regulation National Market System''~\citep{sec2005} to promote ``vigorous competition among markets’’. Today, ``in the United States, stocks listed on Nasdaq or the NYSE can be traded on a myriad of other trading platforms: CBOE's four exchanges, IEX, NYSE Arca, and more. These platforms compete fiercely to attract order flow, and the market share of incumbent exchanges has declined steadily in the United States'' \citep*[Chapter~7]{foucault.al.23}. 

This fragmentation naturally raises the questions of how to model the dynamics of trades and quotes across venues in a systematic and consistent manner, and how to optimize the order routing process in the cross section,\footnote{A related important question is how fees and rebates -- assumed to be fixed in \cite{Maglaras2021}, for example -- are determined by competition between liquidity providers and/or exchanges, cf., e.g., \cite*{pagnotta.philippon.18,baldauf.mollner.21,euch2021} and the references therein.} see \cite[Chapter 7]{foucault.al.23} and the references therein for an overview. \cite*{Maglaras2021} (henceforth \citetalias{Maglaras2021}) have proposed a flexible queueing model that allows to model many of the tradeoffs (e.g., between rebates for placing limit orders and the corresponding waiting times) involved in this ``smart order routing problem'' in a flexible manner.\footnote{Other queueing models with multiple exchanges are studied by~\cite{cont.kukanov.17,kreher.milbradt.23}, for example.} Due to a ``state space collapse'' in the fluid limit of their discrete model,\footnote{Classical references for state-space collapses include~\cite{Bramson:1998aa,Williams:1998aa,Bramson_01,Mandelbaum_04}. Various other scaling limits of limit order books have also been considered. For example, \cite*{Kelly2018,Gao2018} also consider fluid limits, whereas \cite*{cont2013,Horst2019,Hambley2020,Almost2023} study diffusion limits.} the model nevertheless remains remarkably tractable.

However, in the mathematical analysis of this model, two key questions remained open. First, the fluid limit of the system was only motivated informally but not linked to the original queueing system by a rigorous convergence result. Here, the main challenge is that the state dynamics exhibit a singularity when the queue lengths in the system become small. Second, most of the analysis of \citetalias{Maglaras2021} focuses on the steady-state equilibrium of the fluid limit, but stability of the system (i.e., convergence to this stationary value) could only be established for the case of two exchanges. Here, the challenge is that the fluid system is described by a system of coupled nonlinear ODEs. For two exchanges, \citetalias{Maglaras2021} could establish stability using a direct geometric argument, where the state space $\mathbb{R}_{+}^2$ is partitioned into nine regions that can then be dealt with case by case. However, this approach does not generalize when the dimension of the state space increases for a larger number of exchanges. 

In the present study, we first prove that the deterministic fluid system indeed arises as the limit of the discrete queueing model, in the ``fluid scaling’’ of small and frequent orders. To deal with the singularity of the system for short order queues, we first construct an approximating sequence of queueing systems, whose states are bounded away from the origin, and prove that these converge to the conjectured fluid limit. Using a coupling argument, we then show that the original queueing system has the same scaling limit. This convergence theorem also yields wellposedness of the nonlinear ODEs that describe the dynamics of the fluid system.\footnote{Classical wellposedness results such as \cite[Chapter 17.1]{smale2004} only cover such systems if the coefficents are $\mathcal{C}^1$.}

We then turn to the stability of the fluid system and show that the fluid equations for the queue dynamics are \emph{locally} asymptotically stable for an arbitrary number of exchanges. This means that the system converges to its unique stationary point over time if started sufficiently close from this equilibrium value.  The proof of local stability is obtained by analyzing the spectrum of the Jacobian matrix of the states, specifically, by showing that the Jacobian has no negative eigenvalues. In \citetalias{Maglaras2021}, where the proof was carried out for two exchanges, the sign of the real part of the eigenvalues is inferred by the signs of the trace and the $2\times 2$ determinant. This method is difficult to apply for higher dimensions by using minors. In Theorem~\ref{thm:LAS} we therefore instead establish local stability for the system by finding an explicit expression for the eigenvalues of the $N\times N$ Jacobian matrix.

Not surprisingly, \emph{global} stability of the nonlinear fluid system starting from an arbitrary initial configuration is a much more delicate property. Indeed, standard machinery for proving global stability in nonlinear dynamical systems via Lyapunov stability theory (see, e.g., \cite[Chapter 9]{smale2004}) cannot be applied to the system at hand. However, we can report some first progress on this very challenging problem under the additional assumption that the rate at which each exchange attracts market orders only depends on its current queue length but no other idiosyncratic characteristics. Under this condition, we develop a tailor-made geometrical argument, which scales well with the dimensionality of the system, and allows to establish global stability. More specifically, in the proof of Theorem~\ref{thm:GAS} we first show that a close neighborhood of a hyperplane determined by the total mass of the equilibrium point, is an attractor. Then, we prove that once the state enters this neighborhood in the phase space, it must be attracted to the equilibrium point (see Figure~\ref{fig:pr3 aid}). Our additional homogeneity assumption is not required in the two-exchange stability result of \citetalias{Maglaras2021}. It allows us to deal with an arbitrary number of exchanges, while retaining flexibility in choosing heterogeneous values for the rebates they offer. However, it implies that -- in the stationary equilibrium state -- waiting times are uniform across all exchanges. Relaxing this condition therefore is an important but challenging direction for further research.

\paragraph{Structure of the paper:} In
 Section \ref{sec:model} we describe the order-routing model of \citetalias{Maglaras2021}. Section \ref{sec-res} presents our main results regarding the fluid scaling limit of the discrete model and its asymptotic stability. Section~\ref{sec:pr1} is dedicated to the proof of the fluid limit; Section~\ref{sec:proofs} contains the proofs of the stability results.
  
 \paragraph{Notation:} The following notation and conventions are used throughout the paper.
\begin{itemize} 
\item For an integer $N\in\N$ we define $[N]=\{1,2,...,N\}$.
\item The dot product between two vectors $x,y\in\R^N$ is denoted by $x\cdot y=\sum_{i=1}^Nx_iy_i$.
\item The notation $|\cdot|$ is context dependent; when applied to vectors in $\R^N$ it is some norm (since all norms are equivalent it does not matter which one), when applied to sets it is the Lebesgue measure of the set.
\item We use the convention $\inf\emptyset=\infty$ and $\sup\emptyset=-\infty$. 
\item For $a,b\in\R$, we set $a\lor b=\max\{a,b\}$, $a\land b=\min\{a,b\}$, and $a^+=a \lor 0$.
\end{itemize} 
Moreover, we fix a filtered probability space $(\Omega, \mathcal{F},\{\mathcal{F}_t\}_{t\geq0},\PP)$ satisfying the usual conditions of right-continuity and completeness. 
 
\section{The Order-Routing Model of \cite{Maglaras2021}}\label{sec:model}

\subsection{Fragmented Market}

A single risky asset is traded on $N$ separate exchanges labeled by $i=1,\ldots,N$. For simplicity, only the top level on either the bid- or the ask sides of the corresponding limit order books is considered. This means that each venue is modeled by a single FIFO queue $Q^i_t$ describing the number of limit order currently waiting for execution. The vector of all limit order queues is denoted by $Q_t=(Q^i_t)_{i\in[N]}$.

\subsection{Arrivals of Market Orders}

Limit orders on the exchanges are executed when a matching market order is placed. These arrive with an overall rate $\mu>0$ and are routed to the different exchanges with a preferences for more liquidity (i.e., longer queues of limit orders). More specifically, the rate at which market orders arrive at exchange $i$ is 
\be \label{mu-def} 
    \mu_i(Q_t) %=\mu\frac{\beta_i Q^i_t}{\sum_{j=1}^N\beta_j Q^j_t}
    =\mu\frac{\beta_i Q^i_t}{\beta\cdot Q_t}, \quad i\in [N],
\ee
for a given weight vector $\beta= (\beta_1,...,\beta_N) \in\R^N_{+}$. Put differently, the arrival rate on each exchange is given by the fraction of total liquidity available there weighted with the idiosyncratic parameters $\beta_i$, which model that some exchanges may be more attractive than others for characteristics other than liquidity.

Formally, these arrival rates can be realized as follows. Consider the compound Poisson processes
 \begin{align}\label{221}
    M^{i}(t)=\sum_{k=1}^{\wt N^{i}_t}V^i_k, \quad t\geq 0, \, i\in[N],
\end{align}
where the jump times are modeled by unit-rate Poisson processes $\wt{N}^{i}$ and the jump sizes are described by iid $\N$-valued random variables $\{V^i_k\}_{k \geq 1}$, independent of each other and of $\wt{N}^{i}$, with probability mass function $P^i_V$, mean $v$ and finite second moment. The state-dependent arrival rates are in turn incorporated through the time change
\begin{align}\label{222}
    \eta_i(t)=\int_0^t\mu_i(Q_s)ds,\quad t\geq0,\, i\in[N],
\end{align}
i.e., the cumulative number of market orders placed until time $t$ on exchange $i$ is
\be \label{dem} 
    D^i_t = M^i(\eta_i(t)),\quad t\geq0, \, i\in[N]. 
\ee

\subsection{Routing of Limit Orders}

Limit orders are routed to the different exchanges for two reasons. On the one hand, there is a dedicated flow of limit orders that is fully routed to each exchange $i$, modeled by a Poisson process $N^{d,i}$ with intensity $\lambda_i>0$. This part of the order flow represents traders that only provide liquidity on a single exchange (e.g., due to limited infrastructure).

\begin{figure}
    \centering
    \begin{tikzpicture}[font=\footnotesize]
 \foreach \x in {1,2,3,4}{
    \draw[thick] (\x*3,1.5)--(\x*3,0)--(\x*3+1,0)--(\x*3+1,1.5); %buffers
    %\draw[thick] (\x*3,0) rectangle (\x*3+1,1.5);
    \draw[thick] (\x*3,0.3) -- (\x*3+1,0.3); %limit orders
    \draw[thick] (\x*3,0.5) -- (\x*3+1,0.5); %limits orders
    \draw[<-,thick] (\x*3+0.5,1.6) -- (\x*3+0.5,2.1) node[below left] {$\lambda_{\x}$}; %dedicated streams
    \draw[->,thick] (\x*3+0.5,-0.2) -- (\x*3+0.5,-0.7) node[above left] {$\mu_{\x}(Q)$}; %departures
    \draw[->,thick] (\x*0.3+7.25,3.5)--(\x*3+0.6,2.6); %diluted optimized
}
\draw[->,thick] (8,5) -- (8,4.3) node[above left] {$\Lambda$}; %optimized stream
\node at (8,4) {\big[routing to $i^*$ in \eqref{1}\big]};
\draw[->,thick] (10,4)--(11,4) node[right] {$i=0$}; %i=0
%\draw[thick,o-o] (8,3.3)--(8,2.7);
\end{tikzpicture}
    \caption{A system of $N=4$ exchanges where each exchange is presented by a single queue of limit orders. Incoming orders arrive either from dedicated investors at rate $\lambda_i$ or from optimizing investors at rate $\Lambda\chi_i(Q)$. In each venue limit orders are matched with market orders at rate $\mu_i(Q)$.
Exchange $i=0$ corresponds to market orders executed in any of the exchanges.}    \label{fig:model}
 \end{figure}
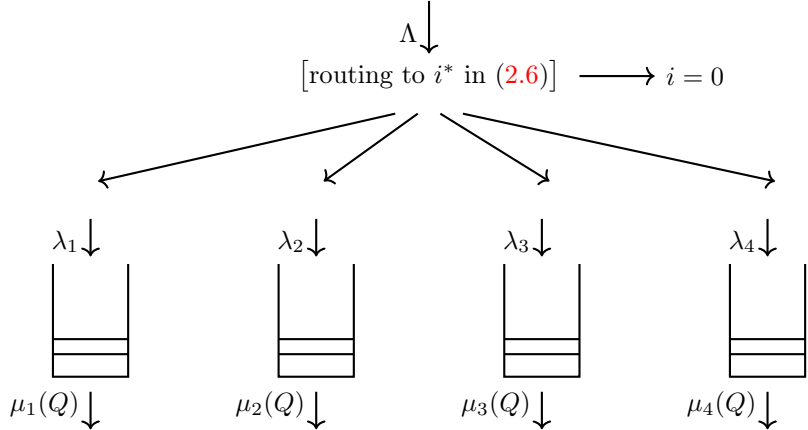

On the other hand, another Poisson process $N_t^o$ with intensity $\Lambda>0$ describes the arrival of limit orders that are \emph{optimally} routed to the exchange that currently offers the best cost-delay tradeoff. Optimized orders can be alternatively sent as market orders and executed immediately if investors find this option favorable. The index $i=0$ in Figure \ref{fig:model} represents this together with the other options. 

The cost-delay tradeoff is modeled in reduced form as follows. In order to encourage liquidity provision, each exchange $i$ offers a \emph{rebate} $r_i \geq 0$ for the placement of limit orders; to avoid having to introduce tiebreakers, $r_i \not = r_j$ for $i \not = j$. Conversely, market orders incur a negative rebate (i.e., a fee) $r_0<0$.

Whereas limit orders are more attractive from a cost perspective, they have the disadvantage of not being executed immediately. This is proxied by the ``expected delay'' in each queue $i$, defined as\footnote{This terminology of \citetalias{Maglaras2021} is exact when the system fluctuates around the stationary state of its fluid limit. Then, the inflows and outflows balance and the time it takes a new order arriving at a queue with a given length and constant service rate indeed has this waiting time.} 
\be 
\sfE\sfD_t^i = \begin{cases}
 \frac{Q^i_{t-}}{v \mu_i(Q_{t-})}\indicator{Q^i_{t-}>0}  &\textrm{for }  i\in[N],   \\
  0 &\textrm{ for }  i=0. 
\end{cases} 
 \ee

The tradeoff between rebates and expected delays is modulated by \emph{type} $\gamma$ of the investors placing the limit orders. This positive variable is determined by iid drawn from an atomless distribution with cumulative distribution function (CDF) $F$. An investor with type $\gamma$ in turn routes their order to the exchange that optimizes the rebate-delay tradeoff\footnote{As $\gamma$ is a continuous random variable and the rebates $\{r_i\}_{i\in [N]}$ are different, there is a unique maximizer with probability one.  } 
 \begin{align}\label{1}
    i^*(\gamma,Q_{t-})\in\argmax_{i\in[N]\cup\{0\}} \left\{\gamma r_i-\sfE\sfD_t^i\right\}.
\end{align}

Denoting by $\tau_k$ the $k$-th jump time of the optimized limit order flow process $N^o$ and by $\gamma_k$ the corresponding trader's type, $i_k^*=i^*(\gamma_k,Q_{\tau_k-})$ is the optimal routing of the limit order arriving at time $\tau_k$. The arrival processes $A^{d,i}$ and $A^{o,i}$ of the dedicated and optimized limit orders at exchange $i$ in turn are
\be \label{arrival} 
\begin{aligned}
    &A^{d,i}_t=\sum_{k=1}^{N^{d,i}_t}B^{d,i}_k,
    &&A^{o,i}_t=\sum_{k=1}^{N_t^o}B^o_k\indicator{i^*_k=i}, \quad i\in [N].
\end{aligned}
\ee
Here, both order sizes $B^{d,i}_k$, $B^{o}_k$ are iid and also independent of the other primitives $V^i_k,N^{d,i}$, $N^o$,$M^i$ of the model. We denote by $P^{d,i}_B$ and $P^o_B$ the probability mass functions of $B^{d,i}$ and $B^o$, and their finite means by $b^{d,i}$, $b^o$, respectively; $B^{d,i}$ and $B^o$ are assumed to have finite second moments as well. 

\begin{remark}
The fluid limit of the discrete system studied in Section~\ref{sec-res} only depends on the average order sizes, whereas fluctuations around these values are averaged out. 
\end{remark}

\subsection{Queue Dynamics}

In summary, the queue length at exchange $i \in [N]$ is 
\begin{align}\label{5}
    Q^i_t=(Q^i_0+A^{d,i}_t+A^{o,i}_t-D^i_t)^+, \quad t \geq 0. 
\end{align}
Here, the application of the operator $(\cdot)^+$ means that market orders that exhaust the queue discard the residue of their order.   

For a vector $q\in\R^N_+$, we define the set $\calS^i(q)$ and the function $\chi_i(q)$:
\begin{align}\label{3}
    &\calS^i(q)=\{\gamma>0:i=i^*(\gamma,q)\}, \qquad \chi_i(q)=\int_{\calS^i(q)}dF(\gamma),  \quad i\in [N].
\end{align}
To wit, $\calS^i(q)$ contains the types of investors that prefer to invest in exchange $i$ when the state of the order queues is $q$, and $\chi_i(q)$ is the probability that an optimized investment would be routed to exchange $i$ when the state is $q$.

\section{Main Results} \label{sec-res}

\subsection{The Fluid Scaling Limit}

In this section we prove that, in the ``fluid scaling'' of small and frequent orders, the discrete queueing system from the previous section indeed converges to the system of nonlinear ODEs suggested in Section 2.3 of \citetalias{Maglaras2021}.

\paragraph{Fluid Scaling}

To this end, we consider a sequence of rescaled systems, indexed by $n\in\N$. The states $\{Q^{i,(n)}\}_{i\in[N]}$ of the $n$-th system have the same dynamics as in the previous section, except that the arrival rates are rescaled with the parameter $n$ to  $n\Lambda$, $n\lambda_i$ and $n\mu_i$ for $i\in[N]$, respectively. 

To counteract this more frequent arrival of new orders, order sizes or, equivalently, the corresponding queue lengths,  are rescaled by $1/n$ to 
\be \label{q-n-def} 
\ol Q^{(n)} = (\ol {Q}^{1,(n)}, ...,\ol {Q}^{N,(n)}),  \quad \textrm{where }  \ \ol {Q}^{i,(n)}=n^{-1}Q^{i,(n)}.
\ee
We note that the functions $\mu_i$ in \eqref{mu-def} and $i^*$, $\chi_i$ in \eqref{3} of the rescaled process can be written as 
\begin{gather}
    \mu_i( {Q}^{(n)}_t)   =\mu\frac{\beta_i\ \ol{Q}^{i,(n)}_{t}}{\beta\cdot \ol{Q}^{(n)}} = \mu_i(\ol{Q}^{(n)}_t), \qquad \sfE\sfD^{i,(n)}_t =\frac{{Q}^{i,(n)}_{t-}}{nv \mu_i( {Q}^{(n)}_{t-})}\indicator{i\neq 0}, \label{f1} \\
    i^*(\gamma,Q^{(n)}_t)\in\argmax_{i\in[N]\cup\{0\}}\{\gamma r_i-\sfE\sfD_t^{i,(n)}\},\qquad \chi_i(Q^{(n)}_t)=\chi_i(\ol{Q}^{(n)}). \notag
\end{gather}

\begin{remark} \label{rem-w}
Recall that $i^*$, $\calS^i$ and $\chi_i$ in \eqref{3} and \eqref{f1} are defined as functions on $\R_+^N$, 
    however from \eqref{f1} it follows that $\sfE\sfD^i$ depends only on $W_t=\beta\cdot Q_t$ and hence can be regarded as a function on $\R_+$.
    Therefore, in the following we often write $i^*(\gamma,W)$, $\calS^i(W)$, and $\chi_i(W)$.
    The process $(W_t)_{t\geq0}$ was introduced in \citetalias{Maglaras2021} as the \emph{workload process}, which fully captures the state of the system in equilibrium (a ``state-space collopase'').
\end{remark}

\paragraph{The Fluid Limit}

\citetalias{Maglaras2021} suggest that as $n \rr \infty$, the rescaled queueing systems $\{\ol{Q}^{n}\}_{n\geq 1}$ convergence to a deterministic limit $\calQ = (\calQ^1,...,\calQ^N)$, described by a system of coupled nonlinear ODEs, where routing probabilities are replaced by appropriate fractions of orders flowing to each exchange: 
\be\label{eq:FE} \calQ_t^i=\calQ^i_0+b^{d,i}\lambda_it+b^o\Lambda\int_0^t\chi_i(\calQ_s)ds-v\int_0^t\mu_i(\calQ_s)ds, \ t\geq 0, \ i\in[N]. 
\ee
Here, $\mu_i$, $\calS^i$ and $\chi_i$ are defined as in \eqref{mu-def} and \eqref{3} and $ i^*$ is given by 
\begin{align}\label{2}
    i^*(\gamma, \calQ)\in\argmax_{i\in[N]\cup\{0\}} \left\{\gamma r_i-\frac{\calQ^i_t}{\mu_i(\calQ_t)v}\indicator{i>0}\right\}.
\end{align}

\paragraph{Assumptions}

We now collect the assumptions required to substantiate the link between the discrete queueing model and the fluid system~\eqref{eq:FE} by a rigorous convergence theorem. To this end, we fix a vector of initial queue lengths $\calQ_0= (\calQ_0^1,...,\calQ_0^N)^\top$ with $\calQ_0^i>0$ for any $i\in [n]$. We also recall the stationary point $\calQ^*= (\calQ^{1*},...,\calQ^{N*})^\top \in \mathbb{R}_{+}^N$ of the fluid system~\eqref{eq:FE} characterized by \citetalias{Maglaras2021} (see their  Eq.~(8)), and define 
$$\calW_0=\sum_{i=1}^N\beta_i\calQ^{i}_0 \quad \mbox{and} \quad \calW^*=\sum_{i=1}^N\beta_i\calQ^{i*}.
$$

By Theorem 2 of \citetalias{Maglaras2021}, we have $\calW^*>0$.
We denote by $\beta_{\min}= \min_{i\in [N]} \beta_i >0$ and $\beta_{\max}= \max_{i\in [N]} \beta_i$ and define the following positive constant:
\be \label{kappa}
\kappa_{\eqref{kappa}} = \beta_{\min}\beta_{\max}^{-1} (\calW_{0} \wedge \calW^*). 
\ee 
Recall that the rebate for market orders is negative ($r_0<0$) and that the rebates $\{r_i\}_{i \in [N]}$ for limit order placements on the different exchanges are nonnegative and different. We set $\calR_i^-=\{j\in[N]:r_j<r_i\}$ and $\calR_i^+=\{j\in[N]:r_j>r_i\}$ and introduce
    \begin{align}\label{898}
        a^+_i=\min_{j\in\calR_i^+}\frac{1}{r_j-r_i}\Big(\frac{1}{\mu\beta_jv}-\frac{1}{\mu\beta_iv}\Big),\qquad
        a^-_i=\max_{j\in\calR_i^-}\frac{1}{r_i-r_j}\Big(\frac{1}{\mu\beta_iv}-\frac{\indicator{j\neq0}}{\mu\beta_jv}\Big),
    \end{align}
as well as
\be \label{amin}
a^-_{\eqref{amin}} =  \min_{i\in [N]}a^-_i>0.
\ee 
%Note that $a^-_{\eqref{amin}}>0$. 

\begin{remark}
The constants $\kappa_{\eqref{kappa}}$ and $a^-_{\eqref{amin}}$ both play a key role in the analysis of the model. Specifically, $\kappa_{\eqref{kappa}}$ is the lower bound on the workload process in the fluid limit (see Lemma \ref{lem-bnd-0}) and the product  $a^-_{\eqref{amin}}\kappa_{\eqref{kappa}}$ gives a bound on the type $\gamma$ below which is always optimal to submit market orders in the fluid limit.  
\end{remark}

The following additional assumptions are needed to derive the fluid limit: %\blue{E: there are many assumptions in this type of queueing models, so I thought that it might be useful to reference all of them together when stating a theorem and when discussing them below. I think that (v) which is commented out $a_i^+\geq a_i^-$ is essential and doesn't appear before. Please note that references to Assumption \ref{ass:}(index), will change now throught the paper.}

\begin{assumption}\label{ass:}
    \begin{enumerate}
        \item[\bf{(i)}]  The distribution function $F$ of investor types $\gamma$ has a density $f=F'$ and $\gamma\mapsto\gamma f(\gamma)$ is strictly decreasing on $[a^-_{\eqref{amin}}\kappa_{\eqref{kappa}},\infty)$.
        \item [\bf{(ii)}]  The arrival rates $\lam_i,\Lambda$, the means $\mu,v$ of market order sizes and rate in \eqref{mu-def}, \eqref{dem}, and the means $b^0,b^{d,i}$ of limit-order sizes in \eqref{arrival} satisfy:
        \begin{align*}
            \sum_{i\in[N]}b^{d,i}\lambda_i<v\mu<\sum_{i\in[N]}b^{d,i}\lambda_i+b^o\Lambda.
        \end{align*}
        %\item [\bf{(iii)}]   The processes $N^{i,d,x}$, $N^{i,o,x}$, $M^{i,x}$, $i\in[N]$ and $x\in\N$, are independent.
        \item  [\bf{(iii)}]  The initial queue lengths in the rescaled systems converge, in that 
        $$
        n^{-1}Q_0^{i,n}\to\calQ_0^i \quad \mbox{a.s.} 
        $$
        \item [\bf{(iv)}]  %We assume that $r_0<0$ and that the rebates $\{r_i\}_{i \in [N]}$ are nonegative and different, i.e. $r_i\neq r_j$ for $i\neq j$ and that 
        $a_i^+\geq a_i^-$ for all $i\in[N]$. 
    \end{enumerate}
\end{assumption}

We now briefly discuss the different parts of Assumption \ref{ass:}:

\begin{remark} 
Assumptions \ref{ass:}(i), (ii) and (iv) are crucial for the proof of the wellposedness of the fluid limit for the queueing system \eqref{q-n-def}. Assumption \ref{ass:}(ii) is equivalent to Assumption 1(ii) in \citetalias{Maglaras2021}. The left inequality in Assumption \ref{ass:}(ii) prevents the system from exploding by having service rate higher than arrival intensity, while the right inequality implies a non-trivial fluid limit due to arrival rate of routed orders. 

Assumptions \ref{ass:}(i) and (iv) are used to show that any solution of \eqref{eq:FE} is bounded away from zero (see Lemma \ref{lem-bnd-0}) and therefore does not approach the singularities of the coefficients $\{\mu^i\}_{i \in [N]}$ at the origin. The existence of the density $f$ in Assumption \ref{ass:}(i) corresponds to Assumption 1(i) in the online supplement of \citetalias{Maglaras2021}. The monotonicity of $\gamma f(\gamma)$ for $ \gamma >a^-_{\eqref{amin}}\kappa_{\eqref{kappa}}$ in Assumptions \ref{ass:}(i) is a relaxed version of Definition 2 in the online supplement of \citetalias{Maglaras2021}, which requires monotonicity of $\gamma f(\gamma)$ for all $\gamma >0$. Many standard distributions satisfy the relaxed version of the assumption, e.g., exponential distributions with density $\lam e^{-\lam\gamma}\mathbf{1}_{\{\gamma>0\}}$ for $\lam>a^-_{\eqref{amin}}\kappa_{\eqref{kappa}}$ or half-normal distributions with density $\sqrt{2/(\pi\sigma^{2})}e^{-\gamma^2/(2\sigma^2)}\mathbf{1}_{\{\gamma>0\}}$ for $\sigma >a^-_{\eqref{amin}}\kappa_{\eqref{kappa}}$. Similar assumptions can be made in order to include the Gamma and the half-Cauchy distributions among others. 

%%\citetalias{Maglaras2021} however failed to provide examples for such distributions. Some of the analytical results in this paper (see e.g. Lemma \ref{lem-bnd-0}) are used in order to relax this assumption. Note that a vast class of distributions satisfy (i) by an appropriate choice of parameters. Some examples are the exponential distribution density $\lam e^{-\lam\gamma}\mathbf{1}_{\{\gamma>0\}}$ where $\lam>a^-_{\eqref{amin}}\kappa_{\eqref{kappa}}$ or half-normal distribution density $\sqrt{2/(\pi\sigma^{2})}e^{-\gamma^2/(2\sigma^2)}\mathbf{1}_{\{\gamma>0\}}$ with $\sigma >a^-_{\eqref{amin}}\kappa_{\eqref{kappa}}$. Similar assumptions can be made in order to include the Gamma and the half-Cauchy distributions among others. 

%he assumption of different rebates corresponds to Assumption 2 in \citetalias{Maglaras2021}. This condition together with 
The assumption on $a_i^-$ and $a_i^+$ in (iv) (which corresponds to Assumption 3 in the online supplement of \citetalias{Maglaras2021}) can be interpreted as the requirement that no exchange is dominated by another. More precisely, if $a_i^+ < a_i^-$ for some $i \in [N]$, then venue $i$ will never be chosen by the decision rule for optimal routing of limit orders in \eqref{1} (see Lemma \ref{lem:chi lipshitz} and  \eqref{g23}, in particular).

%Finally, Assumptions \ref{ass:}(iii) is evidently needed in order to establish a scaling limits. %They were stated in \citetalias{Maglaras2021}, where the convergence to the fluid limit was not established.  
\end{remark}

\paragraph{Convergence Theorem}

We can now state our first main result which proves that, as $n \rr\infty$, the rescaled discrete queueing systems $\ol{Q}^{n}$ converge to a fluid limit described by the coupled nonlinear fluid ODEs~\eqref{eq:FE}.

\begin{theorem}\label{LLN}
Under Assumption \ref{ass:}, as $n\rr \infty$ the sequence $\{\ol{Q}^{n}\}_{n\geq 1}$ converges in probability to a deterministic limit $\calQ = (\calQ^1,...,\calQ^N)$, uniformly on compact subsets of $\mathbb{R}_{+}$. The fluid limit $\calQ$ satisfies the nonlinear ODEs~\eqref{eq:FE}.
%\be\label{eq:FE} \calQ_t^i=\calQ^i_0+b^{d,i}\lambda_it+b^o\Lambda\int_0^t\chi_i(\calQ_s)ds-v\int_0^t\mu_i(\calQ_s)ds, \ t\geq 0, \ i\in[N]. 
%\ee
%Here, $\mu_i$, $\calS^i$ and $\chi_i$ are as in \eqref{mu-def} and \eqref{3} and $ i^*$ is given by 
%\begin{align}\label{2}
%    i^*(\gamma, \calQ)\in\argmax_{i\in[N]\cup\{0\}} \gamma r_i-\frac{\calQ^i_t}{\mu_i(\calQ_t)v}\indicator{i>0}.
%\end{align}
\end{theorem}

The proof of Theorem \ref{LLN} is delegated to Section \ref{sec:pr1} for better readability. Here, we highlight some of the contributions of this result:

\begin{remark} \label{rem-const}
In the fluid limit from Theorem~\ref{LLN}, the parameters $b^{d,i}$, $b^o$ and $v$ can be absorbed in the parameters $\lambda_i$, $\Lambda$ and $\mu$, respectively.
Therefore, there is no loss of generality in letting $b^{d,i}=1$ for all $i\in[N]$, and $b^o=v=1$ when studying the fluid limit itself. The fluid equations are therefore identical to the ones presented in Section~2.3 of \citetalias{Maglaras2021}. 

However, establishing convergence is highly nontrivial. Indeed, standard methods assume that the rates are Lipschitz continuous functions of the state (see, e.g., \cite{kurtz1978,Mandelbaum98}), whereas in the model studied here the functions $\mu_i$ in \eqref{mu-def} are singular at the origin.

This requires a delicate treatment of the system that guarantees that the state remains bounded away from the origin. In order to overcome this difficulty we first construct an approximating sequence of queueing systems which is bounded away from zero, and prove that it converges to \eqref{eq:FE} in Proposition \ref{prop-eps}.  Then we prove that the original system in \eqref{q-n-def} must satisfy the same scaling limit by using a coupling argument. 
\end{remark}

\begin{remark}
The well-posedness of the nonlinear and singular ODE system \eqref{eq:FE}, which was left open in \citetalias{Maglaras2021}, is derived in this paper in two steps. In Proposition~\ref{prop-uniq} we prove that \eqref{eq:FE} admits at most one solution, by using analytic tools. The existence of a solution to \eqref{eq:FE} is then deduced in Proposition~\ref{prop-eps}, by showing that the scaling limit of an approximating queueing system (with respect to \eqref{q-n-def}) satisfies~\eqref{eq:FE}.   
 \end{remark}

%\begin{remark} 
%We can interpret the fluid limit in \eqref{eq:FE} as follows. Consider a system of $N$ stations that are processing fluid mass.
%Station $i$ is associated with a rebate value $r_i$ and its mass-processing rate is $\mu_iv$. 
%Mass is entering the $i$-th station from two sources.
%One of them sends fluid with rate $b^{d,i}\lambda_i$.
%As for the other source, it sends fluid into the system with rate $b^o\Lambda$, and fluid generated by this source has a distribution law $F(\gamma)$, $\gamma>0$.
%This distribution assigns a parameter $\gamma$ to each infinitesimal unit of fluid mass.
%This mass is then distributed in the system between the different stations according to its type, that is, according to its associated parameter $\gamma$.
%Mass associated with $\gamma$ is routed to station $i^*$ which satisfies \eqref{2}; if $i^*$ happens to be $0$ we say that the mass that was routed to station $0$  has exited the system.
%$\calQ^i_t$ in \eqref{eq:FE} is the fluid mass in the $i$-th station at time $t$. 
%\end{remark} 

\subsection{Stabilty of the Fluid Limit}

Our next results deal with the stability of the fluid limit, in the sense of convergence as $t\to\infty$. More specifically, the analysis in \citetalias{Maglaras2021} is centred around the stationary equilibrium of the fluid equations~\eqref{eq:FE}, characterized by the algebraic equations
\be \label{equi} 
    \lambda_i+\Lambda\chi_i(\calQ^*)-\mu_i(\calQ^*)=0, \quad  \ t\geq 0, \ i\in[N]. 
\ee
The existence of a unique solution to \eqref{equi} is established in Theorem 3 of \citetalias{Maglaras2021}; the key question in turn is to study whether the fluid system indeed converges to this steady-state equilibrium in the long run.

\begin{definition}
    An equilibrium $\calQ^*$ is called \emph{locally asymptotically stable} if there exists an $\eps$-neighborhood of $\calQ^*$ such that $\|\calQ_0-\calQ^*\|<\eps$ implies $\calQ_t\to\calQ^*$ as $t\to\infty$ for the solution of \eqref{eq:FE} with initial state $\calQ_0$.
\end{definition}

\begin{definition}
    An equilibrium $\calQ^*$ is called \emph{globally asymptotically stable} if the solution $\calQ_t$ of \eqref{eq:FE} for any initial state $\calQ_0$ satisfies $\calQ_t\to\calQ^*$ as $t\to\infty$. 
\end{definition}

Global asymptotic stability is verified in Theorem 4 in Appendix B.1 of \citetalias{Maglaras2021} for the special case of $N=2$ exchanges. Our second main result proves local asymptotic stability for any $N \geq 1$.

\begin{theorem}\label{thm:LAS}
    Under Assumption~\ref{ass:}, the equilibrium $\calQ^*$ from~\eqref{equi} is locally asymptotically stable for any number $N\geq 1$ of exchanges.
\end{theorem}

Our third major result considers the global asymptotic stability of $\calQ^*$ for any $N \geq 1 $. We show that a sufficient condition is that the parameters $\beta_i$ are the same for all exchanges $i$. This means that the relative attractiveness of the limit order queues for incoming market orders only depends on the relative queue lengths, but no other idiosyncratic characteristics of the exchanges. 

%in the case of similar $\beta_i$'s. Note that in this case the demand rates of the exchanges $\{\mu_i(\cdot)\}_{i\in [N]}$ in \eqref{mu-def} are identical and they depend only on the available liquidity in the exchanges. \blue{YS: I wouldn't say they are identical because there are differences in $r_i$ and $\lambda_i$ and therefore in $Q^i$.}

\begin{theorem}\label{thm:GAS}
Suppose that Assumption~\ref{ass:} holds and $\beta_i=\beta_j$ for all $i,j\in[N]$. Then, the equilibrium $\calQ^*$ is globally asymptotically stable for any $N\geq 1$.  
\end{theorem}
 
The proofs of Theorems \ref{thm:LAS} and \ref{thm:GAS} are provided in Section~\ref{sec:proofs}. Here, we compare them with the corresponding results in \citetalias{Maglaras2021}.

\begin{remark}
  In Lemma B.4.3 of the online supplement of \citetalias{Maglaras2021}, local stability is proved for the special case of $N=2$ exchanges. Our stability result generalizes this to any $N\in\N$.
    For the local stability, both our proof and the proof in \citetalias{Maglaras2021} analyze the spectrum of the Jacobian matrix of the states. 
    In \citetalias{Maglaras2021}, the sign of the real parts of the eigenvalues is inferred by the signs of the trace and the $2\times 2$ determinant, a method that is difficult to extend for higher dimensions by using minors.
    In our proof we therefore instead derive an explicit expression for the eigenvalues of the $N\times N$ Jacobian matrix.
 \end{remark}
 
\begin{remark}
    Global stability is proved in the online supplement of \citetalias{Maglaras2021} (see Lemmas B.4.4 and B.4.5 therein) again for the special case where $N=2$. The proof is based on partitioning the two-dimensional state space into 9 regions and analyzing the state evolution in each of them.
    This method clearly does not scale well for higher-dimensional state spaces.
    In the proof of Theorem~\ref{thm:GAS}, we show that the state space $\R^N_+$ can generally be partitioned and analyzed on three regions when $\beta_i=\beta_j$ for all $i,j\in[N]$. We also note that this case is not covered by \citetalias{Maglaras2021} since they use a coordinate transformation which is not invertible in this case. 
 Thus, our result complements \citetalias{Maglaras2021} also for the case $N=2$.
\end{remark}

\section{Proof of Theorem \ref{LLN}} \label{sec:pr1}

To prepare for the proof of Theorem \ref{LLN}, we first establish the following three auxiliary lemmas. 

\begin{lemma}\label{lem:chi lipshitz}
 For any $\rho>0$ the functions $\{\chi_i\}_{i\in [N]}$ in \eqref{3} are Lipschitz continuous on 
    \be \label{d-rho}
    D_\rho:=\{q\in \mathbb R_+  :   \beta\cdot q>\rho\}.
    \ee
\end{lemma}

\begin{proof}
Recall the definition of $\calS^i$ in \eqref{3} and set
\be \label{gf11}
        \calS^i(q)\Delta\calS^i(q')=(\calS^i(q)\cup\calS^i(q'))\Backslash(\calS^i(q)\cap\calS^i(q'))
\ee
Let $\rho>0$. From Assumption \ref{ass:}(i) it follows that the density $f$ of $\gamma$ must satisfy $\sup_{x>\eps}f(x) <\infty$ for any $\eps>0$. In the following, we will prove that, 
\be \label{asum1} 
|\calS^i(q)\Delta\calS^i(q')|\leq c|q-q'|,\quad \textrm{for all } q, q' \in \mathbb{R}_{+}, 
\ee
for some constant $c>0$. We now show that together with~\eqref{3} this implies that there exists a constant $C>0$ such that 
    \begin{align*}
        \Big|\chi_i(q)-\chi_i(q')\Big| &\leq\int_{\calS^i(q)\Delta\calS^i(q')}dF \\
        &\leq C|q-q'|, \quad \textrm{for all } q, q' \in D_\rho. 
    \end{align*}
Whence, it remains to prove~\eqref{asum1}.  To this end, we need to bound the Lebesgue measure of $\calS^i(q)\Delta\calS^i(q')$.
    From \eqref{f1} (with $n=1$), \eqref{mu-def} and \eqref{3} it follows that for any $i\in[N]$, the scalar $\gamma >0$ is in $\calS^i(q)$ if and only if
    \begin{align*}
        \gamma r_i-\frac{\beta\cdot q}{\mu\beta_iv}\geq\Big(\gamma r_j-\frac{\beta\cdot q}{\mu\beta_jv}\Big)\lor (\gamma r_0),\quad \textrm{for all } j \not =i. 
    \end{align*}
    Equivalently,
\be \label{g1} 
        \gamma(r_i-r_j)\geq\beta\cdot q\Big(\frac{1}{\mu\beta_iv}-\frac{\indicator{j\neq0}}{\mu\beta_jv}\Big),\quad \quad \textrm{for all } j \not =i.
\ee
    Recall that $\calR_i^-=\{j\in[N]:r_j<r_i\}$ and  $\calR_i^+=\{j\in[N]:r_j>r_i\}$ and note that $0\in\calR^-_i$ since $r_0<0$. The inequality~\eqref{g1} therefore is tantamount to 
    \begin{align*}
        \beta\cdot q\min_{j\in\calR_i^+}\frac{1}{r_j-r_i}\Big(\frac{1}{\mu\beta_jv}-\frac{1}{\mu\beta_iv}\Big)\geq\gamma\geq\beta\cdot q\max_{j\in\calR_i^-}\frac{1}{r_i-r_j}\Big(\frac{1}{\mu\beta_iv}-\frac{\indicator{j\neq0}}{\mu\beta_jv}\Big),
    \end{align*}
    and in turn
\be \label{g22}
        a^+_i\beta\cdot q\geq\gamma\geq a^-_{i} \beta\cdot q, 
\ee
    with $a^+_i$ and $a^-_i$ as defined in~\eqref{898}.  
    
Similarly, a number $\gamma>0$ is in $\calS^i(q')$ if and only if
\be \label{g23}
        \beta\cdot q'a^+_i\geq\gamma\geq\beta\cdot q' a^-_{i}.
\ee
From \eqref{gf11}, \eqref{g22}, \eqref{g23} and by the fact that the edges of the interval are linear in the argument $q$ it follows that 
    \begin{align*}
        |\calS^i(q)\Delta\calS^i(q')|\leq c|q-q'|,  \quad \textrm{for all } q, q' \in \R_+, 
    \end{align*}
    for some finite $c>0$. This verifies \eqref{asum1} and thereby completes the proof. 
\end{proof}

\begin{lemma}\label{lem:mu Lipschitz}
    For any $\rho>0$ the functions $\{\mu_i\}_{i\in [N]}$ defined in \eqref{mu-def} are Lipschitz continuous on the set $D_\rho$ defined in \eqref{d-rho}. 
    \end{lemma}
    
\begin{proof}
    The functions $\{\mu_i\}_{i\in [N]}$ are differentiable, with gradient
    \begin{align}\label{40}
        \nabla_q \mu_i(q)=\mu\beta_i\frac{(\beta\cdot q)e_i-q_i\beta}{(\beta\cdot q)^2}=\mu\beta_i\frac{1}{\beta\cdot q}e_i-\mu\beta_iq_i\frac{1}{(\beta\cdot q)^2}\beta.
    \end{align}
    Since $\beta_iq_i\leq\beta\cdot q$ (because $\beta_i, q_i\geq0$) and $\beta\cdot q>\rho$ on $D_\rho$, it follows that
    \begin{align*}
        |\nabla_q \mu_i(q)|\leq\mu\frac{\beta_i}{\rho}|e_i|+\mu\frac{1}{\rho}|\beta|, \quad \textrm{for all }q \in D_\rho.
    \end{align*}
This proves the assertion.
\end{proof}

\begin{lemma} \label{lem-bnd-0} 
For any solution $(\calQ_t)_{t\geq 0}$ of~\eqref{eq:FE}, set $\calW_t = \sum_{i=1}^N \beta_i\calQ_t$ for $t\geq0$. Then under Assumption \ref{ass:} and for $\kappa_{\eqref{kappa}}$ from~\eqref{kappa} we have 
$$
\inf_{t \geq 0} \calW_t >  \kappa_{\eqref{kappa}}. 
$$
\end{lemma}

\begin{proof}
From \eqref{g23} and Assumption \ref{ass:}(iv) it follows that $\calS^i(\calW)=\calW[a_i^-,a_i^+]$, where we recall that $\calW=\beta\cdot\calQ$. Hence from the definition of $\chi_i$ in \eqref{3} we obtain
 \begin{align*}
    \chi_i(\calW)=F(\calW a^+_i)-F(\calW a^-_i).
\end{align*}
Recall that $f =F'$ and in turn
\begin{align}\label{41}
    \frac{d\chi_i(\calW)}{d\calW}=a^+_if(\calW a^+_i)-a^-_if(\calW a^-_i), \quad \textrm{for all } \calW>0.
\end{align}
Together, \eqref{41} and Assumptions \ref{ass:}(i),(iii) yield
\begin{align}\label{554}
    \frac{d\chi_i(\calW)}{d\calW}<0,\quad \textrm{for all } \calW>\kappa_{\eqref{kappa}}, \  i\in[N].
\end{align}
The fluid ODEs~\eqref{eq:FE} can be rewritten as follows:  
\be \label{equ12} 
 \dot{\calQ}^i_t=\Psi_i(\calQ_t), \quad      \Psi_i(\calQ)=\lambda_i+\Lambda\chi_i(\calQ_t)-\mu_i(\calQ_t),\qquad i\in[N].
\ee
Summing over $i$ on both sides of \eqref{equ12} and using \eqref{mu-def} gives
\be \label{tot-m}
\sum_{i=1}^N\dot{\calQ}^i_t = \sum_{i=1}^N\lambda_i+\Lambda(1-\chi_0(\calW_t))- \mu,
\ee
where $\chi_0(\calW)=1-\sum_{i=1}^N\chi_i(\calW)$. 
Note that the equilibrium point $\calQ^*$ satisfies  
\be  \label{g11} 
    \mu=\sum_{i=1}^N\lambda_i+\Lambda(1-\chi_0(\calW^*)).
\ee
We define $\calV_t=\sum_{i=1}^N \calQ^{i}_t$,  $\calV^*=\sum_{i=1}^N\calQ^{i*}$ and $\calW^*=\sum_{i=1}^N\beta_i\calQ^{i*}$. Moreover, we set $\beta_{\min}= \min_{i\in [N]} \beta_i >0$ and $\beta_{\max}= \max_{i\in [N]} \beta_i$ and notice that $\beta_{\min}\calV_t\leq \calW_t \leq \beta_{\max}\calV_t$, and that similar inequalities hold for the equilibrium quantities $\calV^*, \calW^*$. Since $\calW^*>0$ by Theorem 2 of \citetalias{Maglaras2021}, we have $\calV^* \geq \beta^{-1}_{\max}\calW^*>0$. 
Also note that $\chi_0(\calW)$ is monotone increasing when $\calW \geq \kappa_{\eqref{kappa}}$ due to \eqref{554}. 

First assume that $\calW_0 \geq \calW^*$. Then if $\calW_t <\calW^*$ for some $t >0$, \eqref{tot-m}, \eqref{g11} and the monotonicity of $\chi_0$ imply that $\dot{\calV_t}>0$ and hence, by the continuity of the state process, $\calW_t \geq \beta_{\min}  \calV^*$ for all $t\geq 0$. Next consider the case where $\calW_0 \in (0, \calW^*)$. By a similar argument as before, \eqref{tot-m}, \eqref{g11} and the monotonicity of $\chi_0$ imply that $\dot{\calV_t}>0$ as long as $\calW_t <\calW^*$. Hence in this case we have $\calW_t \geq \beta_{\min} \calV_t \geq \beta_{\min}\calV_{0} >0$ for all $t\geq0$. We conclude that in both cases
$$
\calW_t \geq  \beta_{\min}\beta_{\max}^{-1} (\calW_{0} \wedge \calW^*), \quad \textrm{for all } t\geq 0.  
$$
The result thus follows from \eqref{kappa} and since $\calW_{0}>0$ by assumption. 
\end{proof}

Together, Lemmas~\ref{lem:chi lipshitz}--\ref{lem-bnd-0} show that the nonlinear system \eqref{eq:FE} admits at most one solution:

\begin{proposition} \label{prop-uniq}
Assume that $\calQ_0^i\geq 0$ for all $i\in [N]$ and $\sum_{i=1}\beta_i\calQ_0^i>0$. Then, there exists at most one continuous solution to system \eqref{eq:FE}. 
\end{proposition}

 The existence of such solution is deduced in Proposition~\ref{prop-eps} below by showing that the scaling limit of a queueing system satisfies \eqref{eq:FE}. 

To work towards existence of a scaling limit, we prove in the following propostion that the sequence $\{\ol{Q}^{(n)}\}$ defined in \eqref{q-n-def} is \emph{$\calC$-tight}, i.e., it is tight and all limit points are almost surely continuous. In the following we denote the space of c\`adl\`ag functions on $\R^+$ by $\calD$ and equip it with the $J_1$ topology. We write $\calC$ for its subspace of continuous functions and denote the Cartesian product of $N$ such spaces by $\calD^N$ and $\calC^N$. Note that for each $n\geq1$ the process $\ol{Q}^{(n)}$ takes values in $\calD^N$, while the fluid limit $\calQ$ takes values in $\calC^N$. 

 \begin{proposition} \label{prop-tight} 
  The sequence $\{\ol{Q}^{(n)}\}_{n\geq 1}$ is $\calC$-tight.
\end{proposition}

\begin{proof}
Recall the processes $A^{o,i}$, $A^{d,i}$ from~\eqref{arrival} and $D^i$ from~\eqref{dem}. For any $n \geq 1$ we define $A^{o,i,(n)}$, $A^{d,i,(n)}$ as $A^{o,i}$, $A^{d,i}$ with arrival rates $n\Lambda$ and $n\lambda_i$, respectively, and $D^{i,(n)}$ in analogy to $D^i$ but with arrival rate $n\mu_i$, for all $i\in[N]$. Moreover, we set 
\be \label{arrival-n} 
   \ol A^{d,i,(n)}_t=n^{-1}A^{d,i,(n)}_t, \quad  \ol A^{o,i,(n)}_t=n^{-1}A^{o,i,(n)}_t, \quad \ol{D}^{(n)}_t= n^{-1} D^{i,(n)}_t\quad t\geq 0, \ i\in[N]. 
\ee
We also use the notation $\ol{A}_t^{d,(n)}$, $\ol{A}_t^{o,(n)}$ and $\ol{D}_t^{(n)}$ for the $N$-vectors with entries~\eqref{arrival-n}. 

Recalling \eqref{5} and \eqref{q-n-def} we observe that in order to prove that $\{\ol{Q}^{(n)}\}$ is $\calC$-tight we may first prove that $\{\ol{A}^{d,(n)}\}_{n \geq 1}$, $\{\ol{A}^{o,(n)}\}_{n \geq 1}$ and $\{\ol{D}^{(n)}\}_{n \geq 1}$ are $\calC$-tight.
From  Theorem 13.2 in Section 13 of \cite{billingsley99} it follows that we need to bound the tail probabilities of the large jumps and the modulus of continuity for the sequences.
Note that from Equation (12.7) of Section 12 in \cite{billingsley99} it follows that we can use the modulus of continuity $w$ defined in Equation (7.1) of Section 7 therein in order to verify the second condition of Theorem 13.2 of \cite{billingsley99}. Specifically, for an arbitrary stochastic process $\{X_t\}_{t\geq 0}$, the modulus of continuity is given by, 
\be \label{mod} 
w_T(X,\delta) = \sup_{|s-t|<\dl } |X(s)-X(t)|, \quad \textrm{for all } \dl \in (0,1). 
\ee
We start by verifying the conditions for $\ol{A}^{d,(n)}$. Recall that $\ol{A}^{d,i(n)}$ is a non-negative and non-decreasing compound Poisson process with mean jump size $b^{d,i}$ (see \eqref{arrival}). 
Therefore, by \eqref{arrival-n} and Markov inequality we get
\be \label{p1} 
        \PP\Big(\sup_{0\leq t\leq T}|\ol{A}^{d,i(n)}_t|>a\Big)\leq\PP(\ol{A}^{d,i,(n)}_T>a)\leq\frac{\EE[A^{d,i, (n)}_T]}{na}=\frac{ \lambda_i Tb^{d,i}}{ a}. 
\ee
Hence, the right hand side of \eqref{p1} converges to 0 as $a\to\infty$:
\be \label{p12} 
 \lim_{a\rr \infty} \limsup_{n}       \PP\Big(\sup_{0\leq t\leq T}|\ol{A}^{d,i(n)}_t|>a\Big) =0. 
 \ee
This proves condition (i) from Theorem 13.2 of \cite{billingsley99}.

In order to verify Condition (ii) of Theorem 13.2 of \cite{billingsley99} we need to study the modulus of continuity $w_T(\ol{A}^{d,i,(n)},\delta)$ for arbitrary small $\delta >0$. Let $\delta \in (0,1\wedge T)$ and consider a partition of the interval $[0,T]$ into disjoint intervals $\mathbb{T}= \{[t_j,t_{j+1})$, $j=0,...,\lfloor 2T/\delta \rfloor +1\}$ of length at most $\delta/2$ such that $0=t_0<t_1<...<t_{\lfloor 2T/\delta \rfloor +1}=T$. Note that any interval $I \subset [0,T]$ of length at most $\delta$ intersects with at most 3 of the intervals in $\mathbb{T}$.
Therefore, if there is some time interval $I \subset [0,T]$ of length $\delta$ such that the  increment of $\ol{A}^{d,i,(n)}$ is larger than $\eps$, then it must be that for at least one of the intervals in $\mathbb{T}$, the increment of $\ol{A}^{d,i,(n)}$ is larger than $\eps/3$. We therefore have for any $\eps>0$, 
\be \label{inc} 
        \left\{\sup_{0\leq s,t\leq T,|t-s|<\delta}|\ol{A}^{d,i,(n)}_t-\ol{A}^{d,i,(n)}_s|>\eps \right\}\subset \left\{\exists [t_j,t_{j+1})\in \mathbb{T} : | \ol{A}^{d,i,(n)}_{t_{j+1}}-\ol{A}^{d,i,(n)}_{t_j}|>\eps/3 \right\}.
\ee
By \eqref{mod} and \eqref{inc} in conjunction with the union bound, we get
\be\label{383}
    \begin{aligned}
            \PP(w_T(\ol{A}^{d,i,(n)},\delta)>\eps)&=\PP\Big(\sup_{0\leq s,t\leq T,|t-s|<\delta}|\ol{A}^{d,i,(n)}_t-\ol{A}^{d,i,(n)}_s|>\eps\Big) \\
        &\leq\sum_{j=1}^{\lfloor 2T/\delta \rfloor +1}\PP \left (|\ol{A}^{d,i,(n)}_{t_{j+1}}-\ol{A}^{d,i,(n)}_{t_j}|>\eps/3 \right).
    \end{aligned}
    \ee
    Since $A^{d,i,(n)}$ is compound Poisson process with jump size $B_k^{d,i}$ (see \eqref{arrival}) and $t_{j+1}-t_j \leq  \dl/2$, it follows that,
    \begin{align*}
        \PP\left(|\ol{A}^{d,i,(n)}_{t_{j+1}}-\ol{A}^{d,i,(n)}_{t_j}|^2>\eps^2/9 \right)&\leq\frac{\EE[(A^{d,i,(n)}_{t_{j+1}}-A^{d,i,(n)}_{t_j})^2]}{n^2\eps^2/9}\\
        & \leq 9\frac{(n\lambda_i\delta/2)\EE[(B^{d,i,(n)}_1)^2]+(b^{d,i})^2(n\lambda_i\delta/2)^2}{n^2\eps^2}.
    \end{align*}
Together with \eqref{383} we obtain that, for any $\eps>0$, 
\be \label{mod1} 
    \lim_{\dl \rr 0}  \limsup_{n}    \PP(w_T(\ol{A}^{d,i,(n)},\delta)>\eps) = 0. 
\ee
This verifies Condition (ii) from Theorem 13.2 of \cite{billingsley99}. Therefore, the sequence $\{\ol{A}^{d,(n)}\}_{n\geq 1}$ is $\calC$-tight.

We now show that $\{\ol{A}^{o,(n)}\}_{\geq 0}$ in \eqref{arrival-n} is tight using a similar argument. For each $n \geq 1$, we define the compound Poisson process
\be
        X_t^{(n)}=\sum_{k=1}^{N^{o,(n)}_t}B^o_k, \quad t\geq 0, 
\ee
where the jumps $\{B^o_k\}_{k\geq 0}$ and the Poison process $N^{o,(n)}$ are defined in \eqref{arrival}. 
Note that \eqref{p12} and \eqref{mod1} hold for $\{\ol{X}^{(n)}\}$ by the same argument which was used for $\{\ol{A}^{d,(n)}\}_{n \geq 1}$.
Since $0\leq \ol{A}^{o,i,(n)}_t\leq\ol{X}^{(n)}_t$ for all $t\geq0$, $\PP$-a.s., we get that \eqref{p12} holds also for $\{\ol{A}^{o,i,(n)}\}_{n\geq 1}$. 
Moreover, since both $\ol{A}^{o,i,(n)}$ and $\ol{X}^{(n)}$ are monotone non-decreasing in $t$, it follows that $\ol{A}^{o,i,(n)}_t-\ol{A}^{o,i,(n)}_s\leq\ol{X}^{(n)}_t-\ol{X}^{(n)}_t$ for all $0\leq s\leq t$, $\PP$-a.s. Hence, \eqref{mod} implies that for every $\eps>0$, 
    \begin{align*}
   \lim_{\dl \rr 0}  \limsup_{n}        \PP(w_T(\ol{A}^{o,i,(n)},\delta)>\eps)\leq    \lim_{\dl \rr 0}  \limsup_{n} \PP(w_T(\ol{X}^{(n)},\delta)>\eps) =0. 
    \end{align*}
We therefore conclude that $\{\ol{A}^{o,(n)}\}$ is $\calC$-tight.

The processes $M^{i}$ in \eqref{221} are also compound Poisson processes, whose jump sizes have finite second moments. For any $n\geq 1$, define 
$$
\ol{M}^{i,(n)}_t = n^{-1} {M}^{i}_{nt}, \quad t\geq 0, \ n \geq 1. 
$$
Then the sequence $\{\ol{M}^{(n)}\}$ satisfies \eqref{p12} and \eqref{mod1} by the same argument which was used for $\{\ol{A}^{d,(n)}\}_{n \geq 1}$. Recall that instantaneous rates $\mu_i$ in \eqref{mu-def}  are bounded by the constant $\mu$ and that $t\mapsto \eta_i(t)$ in \eqref{222} are non-decreasing for $t\geq 0$.  We therefore get from \eqref{222} and \eqref{f1} for all $0\leq s\leq t $, $i\in [N]$ and $n \geq 1$, 
    \begin{align}
        &\eta_i(t)=\int_0^t\mu_i(\ol{Q}^{(n)}_s)ds\leq \mu t,\label{223}\\
        &|\eta_i(t)-\eta_i(s)|=\int_s^t\mu_i(\ol{Q}^{(n)}_s)ds\leq \mu|t-s|.\label{2224}
    \end{align}
Together, \eqref{dem}, \eqref{arrival-n} and \eqref{223} imply
\be \label{2226}
        \ol{D}^{i,(n)}_t\leq \ol  M^{i, (n)}_{ {\mu}t}, \quad \textrm{for all } t\geq 0, \, \PP-\rm{a.s}
\ee
   As a consequence,
    \begin{align*}
    \lim_{a\rr \infty} \limsup_{n}     \PP\Big(\sup_{0\leq t\leq T}|\ol{D}_t^{i,(n)}|>a\Big)\leq \lim_{a\rr \infty} \limsup_{n} \PP(\ol{M}^{i,(n)}({\mu}T)>a)=0. 
    \end{align*}
 This proves Condition (i) in Theorem 13.2 of \cite{billingsley99} for $\ol{D}^{(n)}$.

From \eqref{2224} it follows that if $|t-s|<\delta$ then $|\eta_i(t)-\eta_i(s)|<  {\mu}\delta$. Since both $\ol D^{i,(n)}$ and $\ol M^{i,(n)}$ are monotone non-decreasing processes, \eqref{dem} and \eqref{arrival-n} give
    \begin{align*}
        \sup_{0\leq s\leq t\leq T, \, |t-s|<\delta}|D_t^{i,(n)}-D^{i,(n)}_s|&=\sup_{0\leq s\leq t\leq T, \, |t-s|<\delta}|M^{i,(n)}(\eta_i (t))-M^{i,(n)}(\eta _i(s))|\\
        &\leq\sup_{0\leq s\leq t\leq   {\mu}T , |t-s|<  {\mu}\delta}|M^{i,(n)}(t)-M^{i,(n)}(s)|.
    \end{align*}
Together with \eqref{mod} this implies
    \begin{align*}
     \lim_{\dl \rr 0}  \limsup_{n}       \PP(w_T(\ol{D}^{i,(n)},\delta)>\eps)\leq  \lim_{\dl \rr 0}  \limsup_{n}    \PP(w_{\mu T}(\ol{M}^{i,(n)},\mu\delta)>\eps) =0. 
    \end{align*}
This proves that Condition (ii) from Theorem 13.2 of \cite{billingsley99} is also satisfied and therefore $\{\ol{D}^{(n)}\}_{n\geq1}$ is $\calC$-tight.

By Prohorov's theorem (see, e.g., Theorem 5.1 of \cite{billingsley99}), since the sequences $\{\ol{A}^{d,(n)}\}_{n\geq 1}$, $\{\ol{A}^{o,(n)}\}_{n\geq 1}$ and $\{\ol{D}^{(n)}\}_{n\geq 1}$ are $\calC$-tight, they converge along subsequences.
    Since the function $(\cdot)^+$ is Lipschitz continuous, it follows from \eqref{5} that the sequence $\{\ol{Q}^{(n)}\}_{n\geq 1}$ also converges along these subsequences. As the space $\calD^N$ is Polish (Theorem 12.2, \cite{billingsley99}), Theorem 5.2 in Section 5 of \cite{billingsley99} implies that $\{\ol{Q}^{(n)}\}_{n\geq1}$ is tight.
    Since the limits of $\{\ol{A}^{d,(n)}\}_{n\geq 1}$, $\{\ol{A}^{o,(n)}\}_{n\geq 1}$ and $\{\ol{D}^{(n)}\}_{n\geq 1}$ are in $\mathcal{C}^N$ and $()^+$ is continuous, the limit of $\{\ol{Q}^{(n)}\}$ is also continuous. In summary, the sequence therefore is $\calC$-tight as asserted.
\end{proof}

The next step in the proof of Theorem \ref{LLN} is to show that any convergent subsequence of $\{\ol{Q}^{(n)}\}$ (which exists by Proposition~\ref {prop-tight} and Prohorov's theorem) indeed satisfies the nonlinear fluid ODEs~\eqref{eq:FE}. Recall that the functions $\mu^i(\cdot)$ in \eqref{mu-def} are discontinuous at the origin (see also \eqref{f1}), which makes it challenging to identify the limit of the sequence 
$$
\ol{D}^{i,(n)}_t=n^{-1}M^i_{t} \left(n \int_0^t\mu_i(\ol{Q}^{(n)}_s)ds \right), \quad t\geq 0, \ n\geq 1. 
$$
In order to overcome this issue we consider truncated versions of $\mu_i(\cdot)$ and $\chi_{i}(\cdot)$ by fixing $\eps \in (0,1)$ and defining for any $q=(q^1,...,q^N)^\top$ in $\mathbb{R}^N$:
\be \label{trunc} 
  \mu_{\eps,i}(q) = \mu\frac{\beta_i q^i}{(\beta\cdot q)\vee \eps},  \qquad  \chi_{\eps,i}(q)=\int_{\calS^i\big((\beta\cdot q)\vee \eps \big)}dF(\gamma),     \quad \ i\in [N].
\ee
(Here, we have used the representation of $\mu_i(\cdot)$ and $\chi_{i}(\cdot)$ from Remark \ref{rem-w} in terms of $W=\beta\cdot q$.) We define $(  Q^{ (n),\eps},   A^{d,(n),\eps},   A^{o,(n),\eps}, {D}^{(n),\eps})$ and $(\ol Q^{ (n),\eps}, \ol A^{d,(n),\eps}, \ol A^{o,(n),\eps},\ol{D}^{(n),\eps})$ in a similar way to \eqref{q-n-def} and \eqref{arrival-n}, only using $ \mu_{\eps,i}$ and $\chi_{\eps,i}$ in  \eqref{trunc} instead of $\mu^i(\cdot)$ and $\chi_{i}$ in \eqref{mu-def} and \eqref{3}. Applying a similar argument as in the proof of Proposition \ref{prop-tight} it follows that for any $\eps \in (0,1)$, the sequence $\{(  Q^{ (n),\eps},   A^{d,(n),\eps},   A^{o,(n)}, {D}^{(n),\eps})\}_{n\geq 1} $ is $\calC$-tight.

In the following proposition we prove that by choosing sufficiently small $\eps$, the limit of any convergent subsequence of $\{\ol{Q}^{(n),\eps}, \}_{n\geq 1} $ satisfies the fluid ODEs~\eqref{eq:FE}. (Recall that we have already established uniqueness for these equations in Proposition~\ref{prop-uniq} above.)

\begin{proposition} \label{prop-eps} 
There exists $\bar \eps>0$ sufficiently small such that for every $\eps \in (0,\bar \eps)$, the limit of any convergent subsequence of $\{\ol{Q}^{(n),\eps}\}_{n\geq 1}$ (in probability, uniformly on compact subsets of $\mathbb{R}_+$) satisfies~\eqref{eq:FE}.  
\end{proposition}

\begin{proof}
Let $\eps >0$, to be specified later. To ease notation, throughout the proof the index $n$ is assumed to belong to an infinite subset of $\N$ where convergence in probability, uniformly on compacts of $\{(\ol{Q}^{(n),\eps},\ol{A}^{d,(n),\eps}$, $\ol{A}^{o,(n),\eps},\ol{D}^{(n),\eps})\}_{n\geq 1}$ to a limiting process $(\calQ^\eps, A^{d,\eps}, A^{o,\eps}, D^{\eps})$ holds.  

We define the following processes on $(\Omega,   {\mathcal F}  ,\{\mathcal{F}_t\}_{t\geq0}, \PP)$: 
    \be \label{marts} 
    \begin{aligned}
        &\ol X^{d,i,(n),\eps}_t= \ol{A}^{d,i,(n),\eps}_t-\lambda_ib^{d,i}t  ,\\
        &\ol X^{o,i,(n),\eps}_t= \ol A^{o,i,(n),\eps}_t-\Lambda b^{o}\int_0^t\chi_{i,\eps}(Q^{(n),\eps}_s)ds,\\
        &\ol Y^{i,(n),\eps}_t=\ol D^{i,(n),\eps}_t-v\int_0^t\mu_{\eps,i}(Q_s^{(n),\eps})ds.
    \end{aligned}
    \ee
Next, we show that Lemma 3.1 from Chapter II.3 in \cite{ikeda1989} implies that $\ol{X}^{d,i,(n),\eps}$, $\ol{X}^{o,i,(n),\eps}$ and $\ol{Y}^{i,(n),\eps}$ are martingales. Note that in \cite{ikeda1989} any compound Poisson process is defined through a multi-dimensional point process, where the first coordinate is the jump intensity and the second is the jump-size.
To connect this to our setting, note that as a continuous-time Markov process $\ol Q^{ (n),\eps}$ can be described through its transition rates.
Let $e_i\in\{0,1\}^N$ be the unit vector of zeros with $1$ in the $i$-th coordinate (i.e. the standard basis of $\R^N$).
The transition rates of $\ol Q^{ (n),\eps}$ are then given by
\begin{align}  
        &q\to q+\frac{x}{n}e_i &&\text{with rate}\qquad n\lambda_iP^{d,i}_B(x)+n\Lambda P^o_B(x)\chi_{\eps,i}(q),\notag\\
        &q\to q-\frac{x}{n}e_i\quad\text{with}\quad q_i-\frac{x}{n}e_i>0 &&\text{with rate}\qquad n\mu_{\eps,i}(q)P_V^i(x),\label{q-n}\\
        &q\to q-\frac{x}{n}e_i\quad\text{with}\quad q_i-\frac{x}{n}e_i=0         &&\text{with rate}\qquad n\mu_{\eps,i}(q)\sum_{x\geq nq_i}P^i_V(x). \notag
\end{align}
Using these intensities, we can in turn define the following point processes:   
             \begin{itemize}
            \item $\textbf{A}^{d,i,(n),\eps}$ on $\R_+^2$  with intensity $n\lambda_i dt P_B^{d,i}(dx)$,
            \item $\textbf{M}^{i,(n),\eps}$ on $\R_+^2\times[0,1]$  with intensity $n\mu dt P_V^{i}(dx) du$,
            \item $\textbf{X}^{(n),\eps}$ on $\R_+^3$ with intensity $n\Lambda dt P_B^o(dx) F(d\gamma)$.
        \end{itemize}
    The processes $A^{d,i,(N),\eps}$, $A^{o,i,(n),\eps}$ and $D^{i,(n),\eps}$ can be defined for any $t\geq 0$: 
    \begin{align*}
        &A^{d,i(n),\eps}_t=\int_0^t\int_0^\infty x\textbf{A}^{d,i(n),\eps}(ds \,dx)  ,\\
        &D_t^{i,(n),\eps}=\int_0^t\int_0^\infty\int_0^1x\indicator{\theta\leq\mu_{i,\eps}(Q^{(n),\eps}_{s-})/\mu}\textbf{M}^{i(n),\eps}(ds \,dx \, d\theta )  ,\\
        &A^{o,i,(n),\eps}_t=\int_0^t\int_0^\infty\int_0^\infty x\indicator{i^*(\gamma,Q^{(n),\eps}_{s-})=i}\textbf{X}^{(n),\eps}(ds \,dx \, d\gamma )  .
    \end{align*}
  From Equation (3.8) in Chapter II.3 of \cite{ikeda1989} it follows that the compensator of $A^{d,i,(n),\eps}$ is $b^{d,i}n\lambda_i t$, the compensator of $D^{i,(n),\eps}$ is $nv\int_0^t\mu_{\eps,i}(\ol{Q}^{(n,\eps)}_s)ds$, and the compensator of $A^{o,i(n),\eps}$ is $n\Lambda b^o\int_0^t\chi_{i,\eps}(\ol{Q}^{(n),\eps}_s)ds$ and therefore $\ol{X}^{d,i,(n),\eps}$, $\ol{X}^{o,i,(n),\eps}$ and $\ol{Y}^{i,(n),\eps}$ are local martingales, and in fact square-integrable martingales because their jump intensities are bounded and the jump sizes are square integrable.  
  
By the Burkholder-Davis-Gundy inequality, these martingales converge to $0$ as $n \rr \infty$ in $L^2$ and hence also in probability, uniformly on compacts, if the expectations of their quadratic variations converge to zero.
 For $\ol{X}^{d,i,(n),\eps}$, we use \eqref{arrival-n}, and the expression for the quadratic variation in Eq. (3.9),  Chapter II.3 of \cite{ikeda1989} and the fact that the jump-sizes are independent from the Poisson process $N^{d,i,(n)}$ to get 
    \begin{align}\label{599}
        \EE[[\ol{X}^{d,i,(n),\eps}]_t]=\frac{1}{n^2}\EE\left [\sum_{k=1}^{N^{d,i,(n)}_t} (B_k^{d,i})^2 \right]=\frac{1}{n^2}\EE[N^{d,i,(n)}_t]\EE[(B_1^{d,i})^2]=\frac{1}{n^2}n\lambda_i t\EE[(B_1^{d,i})^2],
    \end{align}
    which converges to 0 as $n\to\infty$. Therefore by Burkholder-Davis-Gundy inequality it follows that for any $T>0$,
     \be \label{x-d-con}
 \lim_{n\rr \infty} \EE\big[\sup_{t \in [0,T]}(\ol{X}^{d,i,(n),\eps}_t)^2\big]  =  0.
\ee
We handle $\ol{X}^{o,i,(n),\eps}$ similarly, by showing that 
    \begin{align*}
        \EE\left[\big[\ol{X}^{o,i(n),\eps}\big]_t\right]=\frac{1}{n^2}\EE\left[\sum_{k=1}^{N^{o,i,(n)}_t} (B^o_k\indicator{i^*_k=i})^2  \right]\leq\frac{1}{n^2} \EE\left[\sum_{k=1}^{N^{o,i,(n)}_t} (B^o_k)^2 \right].
    \end{align*}
Then, we continue as in \eqref{599} to show that the quadratic variation tends to $0$ as $n\to\infty$, and 
 \be \label{x-o-con}
 \lim_{n\rr \infty} \EE\big[\sup_{t \in [0,T]}(\ol{X}^{o,i,(n),\eps}_t)^2\big]  =  0, \quad \mbox{for any $T>0$.}
\ee
Recall that for any $i \in [N]$, 
$$
\ol{D}^{i,(n),\eps}_t=n^{-1}M^i_{t} \left(n \int_0^t\mu_{\eps,i}(\ol{Q}^{(n),\eps}_s)ds \right), \quad t\geq 0, \ n\geq 1. 
$$
Here the $M^i$'s are the compound Poisson processes defined in \eqref{221}. Hence, to handle $\ol{Y}^{i,(n),\eps}$ in \eqref{marts}, we use similar arguments to \eqref{223}--\eqref{2226} and obtain
\begin{align*}
        \EE[[\ol{Y}^{i,(n),\eps}]_t]=\EE\big[[\ol{M}^{i,(n)}]_{\eta_i(t)}\big]\leq\EE\big[[\ol{M}^{i,(n)}]_{n\mu t}\big], \quad \textrm{for all } t\geq0.
\end{align*}
Since $M^{i,(n)}$ is a unit-rate compound Poisson process with finite second moment jump distribution, it follows by a similar argument as in \eqref{599} that the quadratic variation of $\ol{Y}^{i,(n),\eps}$ converges to $0$ as $n\to\infty$. This proves that for any $T>0$
 \be \label{y-con}
 \lim_{n\rr \infty} \EE\big[\sup_{t \in [0,T]}(\ol{Y}^{i,(n),\eps}_t)^2\big]  =  0.
\ee
Since we have shown in \eqref{x-d-con}--\eqref{y-con} that the martingales in \eqref{marts} converge to $0$, considering \eqref{5} it is left to derive the limits of the compensators. 
Recall that $\mu_i\leq\mu$ and $\chi_i\leq 1$ for all $i\in[N]$, and recall that by Lemmas \ref{lem:chi lipshitz} and \ref{lem:mu Lipschitz}, and \eqref{trunc}, the functions $\mu_{\eps,i}$ and $\chi_{i,\eps}$ are Lipschitz continuous. It follows that there exists a constant $C(T)>0$ such that 
\be \label{comp-con}
      \begin{aligned}
        \sup_{t\in [0,T]} \left| \int_0^t\mu_{\eps,i}(\ol{Q}^{(n),\eps}_s)ds-\int_0^t\mu_{\eps,i}(\calQ_s)ds \right|&\leq C(T)\sup_{t\in [0,T]}|\ol{Q}^{(n),\eps}_t-\calQ_t|  ,\\
        \sup_{t\in [0,T]}\left|  \int_0^t\chi_{i,\eps}(\ol{Q}^{(n),\eps}_s)ds-\int_0^t\chi_{i,\eps} (\calQ_s)ds\right|&\leq C(T)\sup_{t\in [0,T]}|\ol{Q}^{(n),\eps}_t-\calQ_t|.
    \end{aligned}
    \ee
Since $\{\ol{Q}^{(n),\eps}\}_{n\geq 1}$ converge in probability uniformly on compact subsets of $\R_+$, it follows that 
\be\label{comp-con2}
      \begin{aligned}
      \lim_{n\rr \infty}   \int_0^t\mu_{\eps,i}(\ol{Q}^{(n),\eps}_s)ds &= \int_0^t\mu_{\eps,i}(\calQ_s)ds,  \\
          \lim_{n\rr \infty}      \int_0^t\chi_{i,\eps}(\ol{Q}^{(n),\eps}_s)ds &= \int_0^t\chi_{i,\eps}(\calQ_s)ds, 
\end{aligned}
    \ee
where the convergence in \eqref{comp-con2} is also in probability, uniformly on compact subsets of $\R_+$. By writing the rescaled version of \eqref{5} using \eqref{arrival-n} with the truncation $\mu_{\eps,i}$ in \eqref{trunc}, we obtain 
\be \label{res-cont} 
    \ol Q^{i,(n),\eps}_t=(Q^{i,(n)}_0+\ol A^{d,i,(n),\eps}_t+\ol A^{o,i,(n),\eps}_t- \ol D^{i,(n),\eps}_t)^+, \quad t \geq 0. 
\ee
From Assumption \ref{ass:}(iii), \eqref{marts}, \eqref{x-d-con}--\eqref{y-con} and \eqref{comp-con2}, and by using the fact that  $(\cdot)^+$ is Lipschitz continuous, we get that the right-hand side of \eqref{res-cont} converges in probability to the right-hand side of 
\be \label{q1t} 
 \calQ^{i,\eps}_t=       \Big(\calQ^{i}_0+\lambda_ib^{d,i}t+\Lambda b^{o}\int_0^t\chi_{i,\eps}(\calQ^\eps_s)ds-v\int_0^t\mu_{\eps,i}(\calQ^\eps_s)ds\Big)^+.
\ee
Note that by the assumptions of this proposition, $\{\ol{Q}^{(n),\eps}\}_{n\geq 1}$ on the left-hand side of \eqref{res-cont} converges in probability, uniformly on compact subsets of $\mathbb{R}_+$ to $\calQ^{i,\eps}$. 

Now consider the system of ODEs     
\be \label{q0t} 
 \wt{\calQ}_t^{i,\eps} = {\calQ}^{i}_0+\lambda_ib^{d,i}t+\Lambda b^{o}\int_0^t\chi_{i,\eps}(\wt{\calQ}^\eps_s)ds-v\int_0^t\mu_{\eps,i}(\wt{\calQ}^\eps_s)ds , \quad t \geq 0, \ i\in [N]. 
\ee
Note that \eqref{q0t} coincides with the equation~\eqref{q1t} satisfied by $\calQ^{i,\eps}$ up to the first time that $\calQ^{i,\eps}$ hits $0$. However by $\calC$-tightness, $\calQ^{\eps}$ is continuous, and so are the processes inside the $(\cdot)^+$ in \eqref{q1t}. Therefore, the argument of $(\cdot)^+$ in \eqref{q1t} cannot become negative since by \eqref{trunc}, $\mu_{\eps,i}(\calQ^{\eps})=0$ when $\calQ^{i,\eps}=0$. It follows that $\calQ^{i,\eps}$ satisfies \eqref{q0t}.
 
Moreover, from Lemma \ref{lem-bnd-0} and \eqref{trunc} it follows that by choosing $\eps$ sufficiently small, $\calQ^{i,\eps}$ satisfies
    \begin{align*}
        \calQ^{i,\eps}_t= \calQ^{i}_0+\lambda_ib^{d,i}t+\Lambda b^{o}\int_0^t\chi_i(\calQ^{\eps}_s)ds-v\int_0^t\mu_i(\calQ^{\eps}_s)ds , \quad \textrm{for all } t\geq 0, \ i\in [N],
    \end{align*}
which coincides with \eqref{eq:FE}.  By the uniqueness of the solutions to \eqref{eq:FE}, which was proved in Proposition \ref{prop-uniq}, we get the result.  
\end{proof}

We can now complete the proof of the fluid limit results from Theorem \ref{LLN}. 

\begin{proof}[Proof of Theorem \ref{LLN}]
%From Proposition \ref{prop-tight} it follows that in order to prove Theorem \ref{LLN} we need to show that any converging subsequence of $\{(\ol{Q}^{(n)},\ol{A}^{d,(n)}$, $\ol{A}^{o,(n)},\ol{D}^{(n)})\}_{n\geq 1}$ has the same limit $\calQ$ which is the solution to \eqref{eq:FE}. For the sake of readability we assume throughout the proof the index $n$ belongs to an infinite subset of $\N$ where convergence in probability, uniformly on compacts of $\{(\ol{Q}^{(n)},\ol{A}^{d,(n)}$, $\ol{A}^{o,(n)},\ol{D}^{(n)})\}_{n\geq 1}$ to the limiting process $(\calQ, A^{d}, A^{o}, D^{})$ holds. 
Let $\{\calQ_t\}_{t\geq 0}$ be the unique solution to \eqref{eq:FE}. By Lemma~\ref{lem-bnd-0}, there exists $\rho >0 $ such that $\inf_{t\geq0} \beta\cdot \calQ_t > \rho$. Together with Proposition \ref{prop-eps} it follows that there exists $\eps_1 \in (0,(\rho/2)\wedge \bar \eps)$ such that
\be \label{gf1} 
 \lim_{n\rr \infty} \PP\left(\inf_{t\in [0,T]} \beta\cdot \ol{Q}_t^{(n),\eps_1}  > \rho/2   \right) = 1.
\ee
We define the following events:  
$$
\calA^{(n),\rho} = \left\{\inf_{t\in [0,T]} \beta \cdot\ol{Q}_t^{(n)}  > \rho/2  \right\}, \quad n\geq 1. 
$$
Note that $\ol{Q}_t^{(n,)}$ and $\ol{Q}_t^{(n),\eps_1}$ are indistinguishable on $[0,T]$ if $ \inf_{t\in [0,T]} \beta\cdot \ol{Q}_t^{(n)}> \rho/2$. In particular, on the event $\calA^{(n),\rho}$, \eqref{gf1} implies %\red{($\ol{Q}_t^{(n,)}$ and $\ol{Q}_t^{(n),\eps_1}$ are the same as long as $\beta\cdot\ol{Q}^{(n)}_t>\eps_1$, which is like saying that they are the same as long as $\beta\cdot\ol{Q}^{(n),\eps_1}_t>\eps_1$)}
\be \label{gg21}
\lim_{n\rr \infty} \PP(\calA^{(n),\rho}) =1. 
\ee
Let $\eps \in (0,1)$. By \eqref{gg21}, for all $n$ sufficiently large, 
\be \label{fgt1} 
\begin{aligned} 
\PP\left(\sup_{t\in [0,T]} |\ol{Q}_t^{(n)} - \calQ_t |>   \eps \right) &\leq \PP\left(\left\{\sup_{t\in [0,T]} |\ol{Q}_t^{(n)} - \calQ_t | >  \eps  \right\} \cap  \calA^{(n),\rho}  \right) +\PP\big( (\calA^{(n),\rho})^c\big) \\ 
&\leq  \PP\left( \sup_{t\in [0,T]} |\ol{Q}_t^{(n,\eps_1)} - \calQ_t | >  \eps     \right) +\frac{  \eps}{2}\\ 
&\leq   \eps, 
\end{aligned} 
\ee
where we have used Proposition \ref{prop-eps} in the last inequality. Together with \eqref{fgt1}, it follows that $\{\ol{Q}_t^{(n)}\}_{n\geq 1}$ converges to $\calQ_t$ as $n\rr \infty$, uniformly on compact sets. 
 \end{proof}

 \section{Proofs of Theorems \ref{thm:LAS} and \ref{thm:GAS}}\label{sec:proofs}
 \begin{proof} [Proof of Theorem \ref{thm:LAS}]
 
In order to prove local asymptotic stability, we rewrite the system in \eqref{eq:FE} as in \eqref{equ12}. 
From \eqref{equi} and \eqref{equ12} it follows that $\calQ_t^*$ is an equilibrium point. We will compute the Jacobian matrix $J(\calQ)=\{ {\partial}\Psi_i(\calQ)/ {\partial}\calQ^j\}_{i,j \in [N]}$ and prove that all its eigenvalues have negative real parts, hence we establish stability in a neighborhood around the equilibrium (see Theorem in Chapter 8.5, p.175 of \cite{smale2004}).

The Jacobian matrix $J(\calQ)$ is obtained from the derivatives of $\chi$ and $\mu$, which were derived in   \eqref{40} and \eqref{41}:
\be \label{jaco} 
\begin{aligned}
    &J(\calQ)=\\
    &\resizebox{\linewidth}{!}{$\left(\begin{matrix}
        \Lambda\beta_1\frac{d\chi_1(\calW)}{d\calW}-\mu\frac{\beta_1}{\calW}+\mu\frac{\beta_1^2\calQ_1}{\calW^2} &\Lambda\beta_1\frac{d\chi_2(\calW)}{d\calW}+\mu\frac{\beta_1\beta_2\calQ_2}{\calW^2} &\cdots &\Lambda\beta_1\frac{d\chi_N(\calW)}{d\calW}+\mu\frac{\beta_1\beta_N\calQ_N}{\calW^2}\\
        \Lambda\beta_2\frac{d\chi_1(\calW)}{d\calW}+\mu\frac{\beta_2\beta_1\calQ_1}{\calW^2} &\Lambda\beta_2\frac{d\chi_2(\calW)}{d\calW}-\mu\frac{\beta_2}{\calW}+\mu\frac{\beta_2^2\calQ_2}{\calW^2} &\cdots &\Lambda\beta_2\frac{d\chi_N(\calW)}{d\calW}+\mu\frac{\beta_2\beta_N\calQ_N}{\calW^2}\\
        \vdots &\vdots & &\vdots\\
        \Lambda\beta_N\frac{d\chi_1(\calW)}{d\calW}+\mu\frac{\beta_N\beta_1\calQ_1}{\calW^2} &\Lambda\beta_N\frac{d\chi_2(\calW)}{d\calW}+\mu\frac{\beta_N\beta_2\calQ_2}{\calW^2} &\cdots &\Lambda\beta_N\frac{d\chi_N(\calW)}{d\calW}-\mu\frac{\beta_N}{\calW}+\mu\frac{\beta_N^2\calQ_N}{\calW^2}
    \end{matrix}\right).$}
\end{aligned}
\ee
Note that the Jacobian matrix \eqref{jaco} can be decomposed as follows:
\begin{align*}
    J(\calQ)=A-\frac{\mu}{\calW}B. 
\end{align*}
Here, $B=\text{diag}(\beta_1,..., \beta_N)$ and 
\be \label{A-def} 
    A=\beta\Big(\Lambda\frac{d\chi(\calW)}{d\calW}+\frac{\mu}{\calW^2}u\Big)^\top, 
\ee
for $\chi(\calW)=(\chi_1(\calW),...,\chi_N(\calW))^\top$ and $u = (\beta_1\calQ_1,...,\beta_N\calQ_N)^\top$. 
To show that $J$ has no nonnegative eigenvalues we show that the determinant of 
$$
J(\calQ)-\nu I=A-\frac{\mu}{\calW} B-\nu I, \quad \nu \in \mathbb{C},
$$  
is nonzero for any $\nu \geq 0$. Since $A$ is a rank-1 matrix and $-(\mu/\calW)B - \nu I$ is a diagonal matrix we can find an expression for this determinant. 

To this end, first recall that for any $x,y \in \mathbb{R}^N$, a rank-1 matrix $A=xy^\top$ has eigenvalues 0 and $y^\top x$ and, for a diagonal matrix $D$ with eigenvalues $d_i\neq 0$, the determinant of $\det (A-D)$ is given by (see, e.g., \cite[Lemma~1.1]{Ding07})
\be \label{det-eq} 
    \det (A-D)= \det(D)  \det (D^{-1}A-I)=\Big(\prod_{i=1}^Nd_i\Big)\Big(\sum_{i=1}^N\frac{x_iy_i}{d_i}-1\Big)(-1)^{N-1}.
\ee
To use \eqref{det-eq} in our context we set 
\begin{align*}
    x=\beta,\qquad y=\Lambda\frac{d\chi(\calW)}{d\calW}+\frac{\mu}{\calW}u,\qquad D=\frac{\mu}{\calW}B+\nu I,
\end{align*}
and $A=xy^\top$ as in   \eqref{A-def}. It then follows from \eqref{det-eq} that 
\be \label{deq1} 
\begin{aligned} 
  &  \det(J(\calQ)-\nu I) \\
  &=\Big(\prod_{i=1}^N\Big(\frac{\beta_i\mu}{\calW}+\nu\Big)\Big)\Big(\sum_{i=1}^N\frac{\beta^2_i\calQ_i\mu/\calW^2+\Lambda\beta_i (d\chi_i(\calW)/d\calW)}{\beta_i\mu/\calW+\nu}-1\Big)(-1)^{N-1}.
  \end{aligned} 
\ee
To complete the proof we now argue by contradiction, that is, we suppose that there is a value for $\nu$ with non-negative real part for which  \eqref{deq1} vanishes. (As we then have $\nu\neq-\beta_i\mu/\calW$ for any $i\in[N]$ this implies that \eqref{det-eq} is indeed valid because $D$ is invertible.) The determinant in \eqref{deq1} can only be zero if
\be \label{rf1}
    \sum_{i=1}^N\frac{\beta^2_i\calQ_i\mu/\calW^2+\Lambda\beta_i (d\chi_i(\calW)/d\calW)}{\beta_i\mu/\calW+\nu}=1.
\ee
Note that \eqref{rf1} can only hold for real $\nu$, hence this excludes complex eigenvalues.   
Recall that $\calW, \mu, \beta_i >0$. By \eqref{554} and Lemma \ref{lem-bnd-0}  it follows that  
\be \label{deq2} 
    \sum_{i=1}^N\frac{\beta_i (d\chi_i(\calW)/d\calW)}{\beta_i\mu/\calW+\nu}<0,  \quad \textrm{for all } \nu \geq 0. 
\ee
Moreover, as $\calW=\beta\cdot \calQ$,
\be \label{deq3} 
    \sum_{i=1}^N\frac{\beta^2_i\calQ_i\mu/\calW^2}{\beta_i\mu/\calW+\nu}\leq\Big(\sum_{i=1}^N\frac{\beta_i\calQ_i}{\calW}\Big)\max_{i\in[N]}\frac{\beta^i\mu/\calW}{\beta_i\mu/\calW+\nu}\leq1, \quad \textrm{for all } \nu \geq 0. 
\ee
From \eqref{deq1}--\eqref{deq3} we conclude that  $\det(J(\calQ)-\nu I) <0$ for any $\nu\geq0$ so that $\nu \geq 0$ cannot be an eigenvalue of $J$.
In view of the theorem in Chapter 8.5, p.175 of \cite{smale2004} local stability follows.  
\end{proof} 

\begin{proof}[Proof of Theorem \ref{thm:GAS}]
Without loss of generality we can set $\beta_i=1$ for all $i\in[N]$. In this case 
\be \label{w-sepc}
 \calW_t=\beta\cdot \calQ_t = \sum_{i=1}^N\calQ^i_t, 
 \ee
 is simply the total mass in the system.
Note that with identical $\beta_i$'s, the path $t\mapsto \calW_t$ satisfies the following monotonicity conditions:
\begin{align}
    & \dot \calW_t >0 \quad\text{ when }\mu<\sum_{i=1}^N\lambda_i+\Lambda(1-\chi_0(\calW_t)),\label{224}\\
    & \dot \calW_t <0 \quad \text{ when }\mu>\sum_{i=1}^N\lambda_i+\Lambda(1-\chi_0(\calW_t)).\label{225}
\end{align}
That is, if the arrival rate to servers $i =1,...,N$ is greater than the service rate, the total mass increases; if the arrival rate is less than the service rate, the total mass decreases. To see why \eqref{224} and \eqref{225} are true, note that from \eqref{mu-def} and by summing over $i$ in \eqref{equi} the equilibrium point $\calQ^*$ satisfies  
\be \label{mu11} 
    \mu=\sum_{i=1}^N\lambda_i+\Lambda(1-\chi_0(\calW^*)).
\ee
Here, $\chi_0(\calW)=1-\sum_{i=1}^N\chi_i(\calW)$, is increasing in $\calW$ on $[\kappa_{\eqref{kappa}},\infty)$ due to \eqref{554}. We set $\calW^*=\sum_{i=1}^N(\calQ^{*})^i$. 
Therefore, if $\calW_t>\calW^*$ the inequality in \eqref{225} holds and by \eqref{equ12} and \eqref{w-sepc}, $ \dot \calW_t<0$. Similarly if $\calW_t<\calW^*$ then the inequality in \eqref{224} and $\dot \calW_t>0$. This proves that the hyperplane
\be \label{k-def}
    \calK = \Big\{\calQ\in\R_+^N:\sum_{i=1}^N\calQ^i=\calW^* \Big\}   
\ee
is a positively invariant set, which is also attracting.
We also note that the speed in which $\calW_t$ approaches the hyperplane $\calK$ is 
\be \label{speed} 
|\dot{\calW}_t|=\left |\sum_{i=1}^N\lambda_i+\Lambda(1-\chi_0(\calW_t))- \mu\right|
\ee
by \eqref{equ12} and \eqref{w-sepc}.
From \eqref{554} we have $\sum_{i=1}^Nd\chi_i(\calW)/d\calW<0$ which by using $\chi_0(\calW)=1-\sum_{i=1}^N\chi_i(\calW)$, yields $d\chi_0(\calW)/d\calW>0$ for $\calW \in [\kappa_{\eqref{kappa}},\infty)$ (which is where $\calW$ evolves by Lemma \ref{lem-bnd-0}).
This implies that if  $|\calW_t-\calW^*|>\eps$ for any $\eps>0$, then $|\chi_0(\calW_t)-\chi_0(\calW^*)|>\eps'$ for some $\eps'>0$.
Together with \eqref{mu11} it follows that the speed in \eqref{speed} is bounded away from 0 as long as $|\calW_t-\calW^*|>\eps$ for some arbitrary small $\eps>0$, and therefore $\calW_t$ will reach any $\eps$-neighbourhood of $\calW^*$ in finite time.
See Figure \ref{fig:pr3 aid} for an illustration.
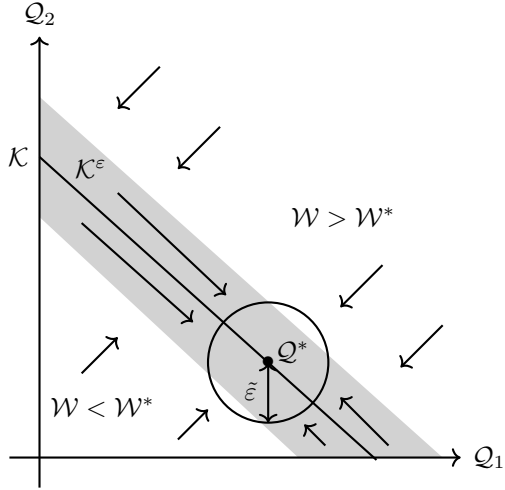
\begin{figure}
    \centering
    \begin{tikzpicture}[scale=0.8,font=\footnotesize]

    % Grey area
  \begin{scope}
    \fill[gray!35] (0,4) -- (0,6) -- (6.7,0) -- (4.3,0) -- (0,4);
  \end{scope}
  % Axes
  \draw[->,thick] (-0.5,0) -- (7,0) node[right] {$\calQ_1$};
  \draw[->,thick] (0,-0.5) -- (0,7) node[above] {$\calQ_2$};
  
  % Function plot
  \draw[thick,domain=0:5.6,smooth,variable=\x] plot ({\x},{5-0.9*\x});
  
  % Annotations
  \node at (0,5) [left] {$\calK$};
  \node at (4.2,1.4) [above] {$\calQ^*$};
  \node at (3.8,1.58) {$\bullet$};
  \node at (1.3,4.8) [left] {$\calK^\eps$};
  \node at (4,4) [right] {$\calW>\calW^*$};
  %circle
  \node at (0,0.5) [above right] {$\calW<\calW^*$};
  \draw[thick] (3.8,1.58) circle (1);
  \draw[<->,thick] (3.8,1.58)--(3.8,0.58) node[left] at (3.8,1.1) {$\tilde \eps$};
  %arrows
  \draw[->,thick] (4.75,0.2)--(4.4,0.55);
  \draw[->,thick] (5.8,0.2)--(5,1);
  \draw[->,thick] (0.7,3.9)--(2.55,2.25);
  \draw[->,thick] (1.3,4.4)--(3.1,2.7);
  
  \draw[->,thick] (0.7,1.4)--(1.3,2);
  \draw[->,thick] (2.3,0.3)--(2.8,0.8);
  \draw[->,thick] (2,6.5)--(1.3,5.8);
  \draw[->,thick] (3,5.5)--(2.3,4.8);
  \draw[->,thick] (5.7,3.2)--(5,2.5);
  \draw[->,thick] (6.7,2.2)--(6,1.5);
\end{tikzpicture}
    \caption{An illustration of the state dynamics on different regions, indicated by the arrows.
    The black diagonal line is the hyperplane $\calK$, the grey area around it is $\calK^\eps$.
    The ball of radius $\tilde \eps$ around the equilibrium point $\calQ^*$ is where local stability holds by Theorem \ref{thm:LAS}.
    On the regions $\calW<\calW^*$ and $\calW>\calW^*$ the state approaches $\calK^{\eps}$. }
    \label{fig:pr3 aid}
\end{figure}

On the hyperplane $\calK$ in \eqref{k-def} and very close to it, the trajectories of the system \eqref{equ12}  are moving towards the equilibrium point $\calQ^*$.
To make this statement precise we define the set, 
\be  \label{k-eps} 
    \calK^\eps=\Big\{\calQ\in\R^N_+: \big \|\sum_{i=1}^N\calQ-\calW^* \big\|\leq\eps\Big\},\qquad\eps>0.
\ee
Now, suppose that the state $\calQ_t$ is in $\calK^\eps$ for some small enough $\eps>0$ to be specified later. 
We define 
\be\label{783}
\begin{aligned}
    \delta_1(\calW_t) &= \Lambda\chi_i(\calW_t)-\Lambda\chi_i(\calW^*),\\
    \delta_2(\calW_t)&= \mu\left(\frac{1}{\calW_t}-\frac{1}{\calW^*} \right). 
\end{aligned}
\ee
From Lemma \ref{lem-bnd-0} we have $\inf_{t \geq 0}\calW_t>\kappa_{\eqref{kappa}}>0$. 
Together with the fact that $\chi_i$ is Lipschitz continuous on $D_{\kappa_{\eqref{kappa}}}$ by Lemma \ref{lem:chi lipshitz}, it follows that there exists $C_2>0$ such that, 
\be \label{f311} 
 |\delta_1(\calW)|\lor|\delta_2(\calW)|<    C_2 | \calW-\calW^*|      , \quad \textrm{for all } \calW \geq \kappa_{\eqref{kappa}}. 
 \ee
Suppose further that $\calQ^i_t >\calQ^{i*}$ for some $i$.
Then from the definition of the equilibrium \eqref{equi},  \eqref{equ12} and \eqref{783} it follows that for all $\calQ_t$ satisfying $|\calW_t- \cal W^*|\leq \cal W^*$ we have,
\be\label{409}
\begin{aligned}
    \dot{\calQ}^i_t&=\lambda_i+\Lambda\chi_i(\calW_t)-\mu\frac{\calQ^i_t}{\calW_t} \\
    &=\lambda_i+\Lambda\chi_i(\calW^*)+\delta_1(\calW_t)-\mu\frac{\calQ^i_t}{\calW^*}-\delta_2(\calW_t) \calQ^{i}_t \\
    &=\Big(\frac{\mu}{\calW^*}+\delta_2(\calW_t)\Big)(\calQ^{i*}-\calQ^i_t)+\delta_1(\calW_t)-\calQ^{i*}\delta_2(\calW_t) \\ 
        &=\frac{\mu}{\calW_t}(\calQ^{i*}-\calQ^i_t)+\delta_1(\calW_t)-\calQ^{i*}\delta_2(\calW_t) \\ 
        &\leq \frac{\mu}{2\calW^*}(\calQ^{i*}-\calQ^i_t)+C_2(1+\calQ^{i*})   | \calW_t-\calW^*| .
\end{aligned}
\ee
Let $\tilde \eps >0$ be sufficiently small such that the local stability Theorem \ref{thm:LAS} holds on the $\tilde \eps$-ball centred at $\calQ^*$, denoted by $ B_{\tilde \eps}(\calQ^*)$. Next, we choose $\eps \in (0,(\tilde \eps/\sqrt{N})\wedge \cal W^*)$ sufficiently small which satisfies the following property: if 
$\calQ^i_t-\calQ^{i*}>\tilde\eps/\sqrt{N}$ and $\calQ_t \in \calK^\eps$ then we have 
 \begin{align*}
    \dot{\calQ^i_t}  &\leq  - \frac{\mu}{2\calW^*} \frac{\tilde \eps}{\sqrt{N} }+ C_2(1+\calQ^{i*})  \eps < 0, \quad \textrm{for all } \calQ_t \in \calK^\eps  \textrm{ with } \calQ^i_t-\calQ^{i*} > \frac{\tilde\eps}{\sqrt{N}}. 
\end{align*}
Such an $\eps$ can indeed be chosen thanks to \eqref{k-eps},  \eqref{409} and by recalling that $ \calW_t = \sum_{i=1}^N \calQ^i_t$.
It follows that the speed $ \dot{\calQ^i_t}$ is bounded away from zero for any $\calQ^i_t>\calQ^{i*} +\tilde\eps/\sqrt{N}$ in $\calK^\eps$. Hence $\calQ^i$ will hit the level $\calQ^{i*}+\tilde \eps/\sqrt{N}$ in a finite time, and will stay within this distance from $\calQ^{i*}$ permanently. The same argument applies for the case where $\calQ^i_t<\calQ^{i*} -\tilde\eps/\sqrt{N}$. We therefore conclude that $\calQ_t$ will enter the local stability region $ B_{\tilde \eps}(\calQ^*)$ in a finite time. The local asymptotic stability property on $ B_{\tilde \eps}(\calQ^*)$, which was obtained in Theorem \ref{thm:LAS}, yields convergence to the equilibrium $\calQ_t\to\calQ^{*}$ as $t\to \infty$. 
\end{proof}

\bibliographystyle{abbrvnat}
\printindex
\bibliography{Bib}

 \end{document}